\documentclass[11pt]{amsart}
\usepackage{amsfonts}
\usepackage{amssymb,latexsym,amsmath,amsthm, mathrsfs}
\usepackage[all]{xy}
\usepackage{bbm}
\usepackage{color}

\def\mycaption#1#2{
    \addtocounter{figure}{1}
    \parbox[b]{#1}{{\scshape Figure $\arabic{figure}$.}\ #2}
}   % used for multiple captions in the same figure environment

\def\pd<#1>{{\!\left<#1\right>}}

\newtheorem{thm}{Theorem}[section]
\newtheorem{lem}[thm]{Lemma}
\newtheorem{prop}[thm]{Proposition}
\newtheorem{cor}[thm]{Corollary}
\newtheorem{df}[thm]{Definition}
\newtheorem{fact}[thm]{Fact}

\newtheorem*{question}{Question}
\newtheorem*{mainthm}{Main Theorem}

\def\ss{\section}
\def\sss{\subsection}

\def\P{{\mathbb P}}
\def\Q{{\mathbb Q}}

\def\O{{\mathcal O}}
\def\Y{{\widetilde{Y}}}
\def\y{{\tilde{y}}}

\def\ave{{\textrm{ave}}}

\def\tm{{\times}}
\def\ts{{\otimes}}%TenSor.
\def\btm{,}
\def\rmk{\medskip\noindent \textbf{Remark: }}

\def\egs{\medskip\noindent \textbf{Examples: }}

\subjclass[2000]{14C15, 14N20}
\title[Chow motive]{Chow Motive of Fulton-MacPherson~configuration~spaces and~wonderful~compactifications}
\author{Li Li}

\begin{document}
\maketitle

\iffalse
\begin{abstract}
\indent We study the Chow groups and the Chow motives of the wonderful compactification
$Y_{\mathcal{G}}$ of an arrangement of subvarieties. We prove a natural decomposition of the Chow
motive of $Y_\mathcal{G}$, in particular of the Fulton-MacPherson configuration space $X[n]$. As a
consequence, we prove a decomposition of the Chow motive of $X[n]/\frak{S}_n$. A generating
function for the Chow groups and for the Chow motive of $X[n]$ is given.
\end{abstract}
\tableofcontents
\fi

\ss {Introduction}

The purpose of this article is to study the Chow groups and Chow motives of the so-called wonderful compactifications of an arrangement of
subvarieties, in particular the Fulton-MacPherson configuration spaces.

All the varieties in the paper are over an algebraically closed field. Let $Y$ be a nonsingular quasi-projective variety. Let $\mathcal{S}$ be an
arrangement of subvarieties of $Y$ (cf. Definition \ref{df:arrangement}). Let $\mathcal{G}$ be a building set of $\mathcal{S}$, i.e., a finite set of
nonsingular subvarieties in $\mathcal{S}$ satisfying Definition \ref{df:building set}. The wonderful compactification $Y_\mathcal{G}$ is constructed
by blowing up $Y$ along subvarieties in $\mathcal{G}$ successively (cf. Definition \ref{df:wonderful compactification}). There are different
orders in which the blow-ups can be carried out, for example we can blow up along the centers in any order that is compatible with the inclusion
relation.  There are many important examples of such compactifications: De Concini and Procesi's wonderful model of a subspace arrangement, the
Fulton-MacPherson configuration spaces, the moduli space $\overline{\mathcal{M}}_{0,n}$ of stable rational curves with $n$ marked points, etc. These
spaces have many properties in common. Studying them by a uniform method gives us better understanding of these spaces. In this article, we study
their Chow groups and Chow motives.

If we assume that $Y$ is projective, then the Chow motive of $Y_\mathcal{G}$, denoted by $h(Y_\mathcal{G})$, can be decomposed canonically into a
direct sum of the motive of $Y$ and the twisted motives of the subvarieties in the arrangement (cf. \S2.1 for a review of Chow motives). We will
prove the following theorem, where the precise definition of the set $M_\mathcal{T}$ and the subvarieties $Y_0\mathcal{T}$ of $Y$ are in \S3.
\begin{mainthm}[Theorems \ref{main thm wonderful}, \ref{main thm correspondence}] Let $Y$ be a nonsingular quasi-projective variety, $\mathcal{G}$ be
a building set and $Y_\mathcal{G}$ be the wonderful compactification $Y_\mathcal{G}$. Then we have the Chow group decomposition
$$A^*Y_\mathcal{G}=A^*Y\oplus\bigoplus_{\mathcal{T}}\bigoplus_{\underline{\mu}\in
M_\mathcal{T}} A^{*-\|\underline{\mu}\|}(Y_0\mathcal{T})$$ where $\mathcal{T}$ runs through all $\mathcal{G}$-nests. Moreover, when $Y$ is projective
we have a Chow motive decomposition
$$h(Y_\mathcal{G})\cong h(Y)\oplus\bigoplus_{\mathcal{T}}\bigoplus_{\underline{\mu}\in
M_\mathcal{T}} h(Y_0\mathcal{T})(||\underline{\mu}||)$$ where $\mathcal{T}$ runs through all $\mathcal{G}$-nests. In this case the correspondences
giving the above isomorphism are canonical in the following sense: although there is no canonical order of blow-ups (in general) to construct
$Y_\mathcal{G}$, the correspondences turn out to be independent of the order we choose.

\end{mainthm}

The Fulton-MacPherson configuration space $X[n]$ is one of the most interesting examples of the wonderful compactification $Y_\mathcal{G}$ where $Y=X^n$ and $\mathcal{G}$ is the the set of all the diagonals in $X^n$ (cf. \S4.1). Applying the main theorem to $X[n]$, we obtain the following
theorem, where the precise definition of the nests $\mathcal{S}$, the polydiagonals $\Delta_\mathcal{S}$, the integers $c(\mathcal{S})$, the sets of lattice points $M_{\mathcal{S}}$, and the correspondences
$\alpha_{\mathcal{S}, \underline{\mu}}$ and $\beta_{\mathcal{S}, \underline{\mu}}$ are in \S4.1.

\medskip\noindent\textbf{Theorem \ref{main FM motive}.}\emph{ Let $X$ be a nonsingular projective variety. Then there is a canonical isomorphism of Chow motives
$$\bigoplus_{\mathcal{S}}\bigoplus_{\underline{\mu}\in M_{\mathcal{S}}}\alpha_{\mathcal{S}, \underline{\mu}}:\quad
h(X[n])\cong \bigoplus_{\mathcal{S}}\bigoplus_{\underline{\mu}\in M_{\mathcal{S}}}h(\Delta_\mathcal{S})(\|\underline{\mu}\|)$$ with the inverse
$\sum\limits_{\mathcal{S}}\sum\limits_{\underline{\mu}\in\mathcal{S}}\beta_{\mathcal{S}, \underline{\mu}}$. Equivalently, we have the following
decomposition of the Chow motive of $X[n]$:
$$h(X[n])\cong \bigoplus_{\mathcal{S}}\bigoplus_{\underline{\mu}\in M_{\mathcal{S}}} h(X^{c(\mathcal{S})})(\|\underline{\mu}\|).$$}

Here are two consequences of this theorem. One is that we can easily express the decomposition of $h(X[n])$ using a generating function $N(x,t)$, as
follows.

\medskip
\noindent\textbf{Theorem \ref{generating}.} \emph{ Define $f_i(x)$ to be the polynomials whose
exponential generating function $N(x,t)=\sum\limits_{i\ge 1}f_i(x)\frac{t^i}{i!}$ satisfies the
identity
$$(1-x)x^dt+(1-x^{d+1})=\exp(x^dN)-x^{d+1}\exp(N)$$
where $d=\dim X$. Then
$$h(X[n])=\bigoplus_{\substack{1\le k\le n\\i\ge 0}}\big{(} h(X^k)(i)\big{)}^{\oplus
[\frac{x^it^n}{n!}]\frac{N^k}{k!}}.
$$}

The other consequence is a decomposition of the Chow motive of the quotient variety $X[n]/{\frak{S}_n}$ obtained from the natural symmetric group $\frak{S}_n$
action on $X[n]$. To make sense of the motive of a quotient variety, we assume the base field is of characteristic 0. The correspondences appeared in
Theorem \ref{main FM motive} are canonical, and therefore symmetric with respect to the symmetric group $\frak{S}_n$. It is then possible to compute
the $\mathfrak{S}_n$-invariant part of $h(X[n])$, which is the Chow motive of $X[n]/\frak{S}_n$. As pointed out by \cite{FM}, unlike the isotropy groups of a point in $X^n$, the isotropy group of
any point in $X[n]$ is always solvable, therefore the singularity of $X[n]/\mathfrak{S}_n$ is ``better'' than the singularity of the symmetric product
$X^{(n)}:=X^n/\frak{S}_n$. It would be interesting to see how different is the Chow motive $h(X[n]/\mathfrak{S}_n)$ from $h(X^{(n)})$. In the
following theorem, an unlabeled weighted forest is a forest whose nodes are not labeled and that each non-leaf node is attached by a positive integer called
weight; we call an unlabeled weighted forest \emph{of type $\nu:=\{n_1,\dots, n_r\}$} if the forest is of the form $n_1T_1+\cdots +n_rT_r$, where $T_i$
are mutually distinct unlabeled weighted tree.

\smallskip \noindent\textbf{Theorem \ref{motive X[n]/Sn}.} \emph{ For any unordered set of positive integers $\nu=\{n_1,\dots, n_r\}$ and any non-negative integer $m$, let
$\lambda(\nu,m)$ to be the number of unlabeled weighted forest with $n$ leaves, of type $\nu$ and of total weight $m$, such that at each non-leaf $v$
with $c_v$ children, the weight $m_v$ satisfies $1\le m_v\le (c_v-1)\dim X-1$. Then
$$h(X[n]/\frak{S}_n)=\bigoplus_{\nu,m}\Big{[} h\big{(}X^{(n_1)}\times \cdots\times X^{(n_r)}\big{)}(m)\Big{]}^{\oplus\lambda(\nu,m)}.$$
}

The importance of all the above results of Chow motives can be seen through a working principle:

\smallskip

 \noindent{\textbf{Principle: }} \emph{A result proved for Chow motives is valid if we replace them
by homological/numerical motives, Chow groups $A^*_\mathbb{Q}$, cohomology groups $H^*_\mathbb{Q}$,
Grothendieck groups (the aforementioned groups are taken with $\mathbb{Q}$-coefficients), Hodge
structures, etc.}

\smallskip

Thus for example, we have a decomposition for the $\mathbb{Q}$-coefficient singular cohomology of $Y_\mathcal{G}$, $X[n]$ and $X[n]/\frak{S}_n$.

\smallskip

The paper is organized as follows. \S2 contains a review of motives and the wonderful compactifications of arrangement of subvarieties. In \S3 a
motivic decomposition for the wonderful compactifications is proved. In \S4 we give a motivic decomposition for the Fulton-MacPherson configuration
spaces. \S5 gives a motivic decomposition for the quotient variety $X[n]/\frak{S}_n$.

\medskip \noindent\textbf{Acknowledgements.} The paper is based on part of the author's thesis. In many ways I am greatly indebted to Mark de
Cataldo, my Ph.D. advisor. I would also thank Blaine Lawson, Sorin Popescu, Dror Varolin, Byungheup Jun for encouragement and useful discussions. The
author thanks the referee for the detailed review of the first version with many very helpful suggestions to clarify and simplify the paper.
%%%%%%%%%%%%%%%%%%%%%%%%%%%%%%%%%%%%%%%%%%%%%%%%%%%%%%%%%%%%%%%%%%%%%%%%%

\ss {Preliminaries}

\sss{Motives}

Given an algebraic variety $X$ of dimension $d$,  let $A^{i}X=A_{d-i}X$ be the Chow group of
codimension $i$, i.e., the group of algebraic cycles of codimension $i$ in $X$ modulo rational
equivalence. Define $A^i_{\Q}X=A^iX\otimes_\mathbb{Z}\Q$.

Let $X, Y$ be two non-singular projective varieties. The group of {\em correspondences of degree
$r$} from $X$ to $Y$ is defined as
$$Corr^r(X,Y):=A^{\dim X+r}(X\times Y).$$
The group $Corr^r_\mathbb{Q}(X,Y)$ denotes the tensor of $Corr^r(X,Y)$ with $\mathbb{Q}$.

The \emph{composition} of two correspondences $f\in Corr^r(X_1,X_2)$ and $g\in Corr^s(X_2, X_3)$ is a correspondence in $Corr^{r+s}(X_1, X_3)$
defined as
$$g\circ f:=\pi_{13*}^{}(\pi_{12}^*f\cdot\pi_{23}^*g)$$
where $\pi_{ij}$ is the projection from $X_1\times X_2\times X_3$ to $X_i\times X_j$.

A correspondence $p\in Corr^0(X,X)$ is called a \emph{projector} of $X$ if $p^2(:=p\circ p)=p$.

  Let $\mathcal{V}$ denote the category of (not necessarily connected) non-singular projective varieties
over a field $k$.

\begin{df}[\cite{CH}] The \emph{category of Chow motives over $k$}, denoted by $CH\mathcal{M}$, is
defined as follows:
an object of $CH\mathcal{M}$, called a {\em Chow motive}, is a triple $(X,p,r)$, where $X$ is a nonsingular projective
variety, $p$ is a projector of $X$ and $r$ is an integer. The morphisms in $CH\mathcal{M}$ are defined as
$$Hom_{CH\mathcal{M}}\big{(}(X,p,r), (Y,q,s)\big{)}:=q\circ Corr^{s-r}(X, Y)\circ p.$$
The composition of morphisms is defined as the composition of correspondences.
\end{df}
For a Chow motive $M=(X,p,r)$ and an integer $\ell$, we define $$M(\ell):=(X,p,r+\ell).$$

There is a natural contravariant functor $h$ from $\mathcal{V}$ to $CH\mathcal{M}$, which sends $X$
to $(X, id_X, 0)$ and sends a morphism $f: X\to Y$ to $\Gamma_f^t: h(Y)\to h(X)$, the transpose of
the graph of $f$. Naturally, $h(X)(\ell)$ stands for the Chow motive $(X,id_X,\ell)$.

According to \cite{dBV}, we can generalize the theory of Chow motives on nonsingular projective varieties to the one on varieties which are quotients
of smooth projective varieties by finite group actions. To be more precisely, let $\mathcal{V}'$ be the category of varieties of type $X/G$ with
$X\in Ob \mathcal{V}$ and $G$ a finite group. We can define the group of correspondences $Corr_\Q^r(X', Y')$ for $X', Y'\in \mathcal{V}'$  and the
category of Chow motives $CH\mathcal{M}'$ similar to the nonsingular case. (The difference is that we have to use $\mathbb{Q}$-coefficients). There
is a natural contravariant functor $h: \mathcal{V}'\to CH\mathcal{M}'$.

Define the $G$-average correspondence $\ave_G$ as $$\ave_G:=\dfrac{1}{|G|}\sum [g]\in Corr^0_\Q(X,X)$$ where $[g]$ is given by the graph of $g$ in $X\times X$. By \cite{dBV}
Proposition 1.2, there is an isomorphism
 $$h(X/G)\cong (X, \ave\Delta)\cong h(X)^G.$$
Such a definition is consistent with the realization functors and $\Q$-coefficient Chow groups.

\sss{Wonderful compactification of an arrangement of subvarieties}

The wonderful compactification of an arrangement of subvarieties is introduced in \cite{Li} as a
generalization of De Concini and Procesi's wonderful model of subspace arrangements. We briefly
review the definition and some properties of such compactifications. For details we refer to
\cite{Li}.

\begin{df}\label{df:arrangement} A (simple) {\em arrangement} of subvarieties of\, $Y$ is a finite set $\mathcal{S}=\{S_i\}$ of
nonsingular closed subvarieties of\, $Y$ satisfying the following conditions:
\begin{itemize}
\item[(1)] $S_i$ and $S_j$ intersect cleanly (i.e. their intersection is nonsingular and $T(S_i\cap S_j)=T(S_i)|_{(S_i\cap S_j)}\cap T(S_j)|_{(S_i\cap S_j)}$),
\item[(2)] $S_i\cap S_j$ is either empty or equal to some $S_k\in\mathcal{S}$.
\end{itemize}
\end{df}

\begin{df}\label{df:building set} Let $\mathcal{S}$ be an arrangement of subvarieties of $\, Y$. A subset $\mathcal{G}\subseteq
\mathcal{S}$ is called a \emph{building set of $\mathcal{S}$} if $\forall S\in \mathcal{S}$, the minimal elements in $\mathcal{G}$ which contains $S$
intersect transversally and their intersection is $S$ (this condition is always satisfied if $S\in \mathcal{G}$). These minimal elements are called
the \emph{$\mathcal{G}$-factors} of $S$. We call a finite set $\mathcal{G}$ of subvarieties a \emph{building set} if the set
$$\mathcal{S}:=\{\bigcap_{V\in \mathcal{F}}V\}_\mathcal{F},$$ where $\mathcal{F}$ runs through all subsets of $\mathcal{G}$, is
  an arrangement and $\mathcal{G}$ is a building set of $\mathcal{S}$ (for $\mathcal{F}=\emptyset$ we set $\bigcap_{V\in \mathcal{F}}V=Y$).
   In this case we call $\mathcal{S}$ the \emph{induced arrangement} of $\mathcal{G}$.
\end{df}

\begin{df}\label{nest}
Let $\mathcal{G}$ be a building set. A subset $\mathcal{T}\subseteq\mathcal{G}$ is called \emph{$\mathcal{G}$-nested} (or a
\emph{$\mathcal{G}$-nest}) if it satisfies one of the following equivalent relations:
\begin{enumerate}
\item There is a flag of elements in $\mathcal{S}$: $S_1\subseteq S_2\subseteq\dots\subseteq S_k$, such that
 $$\mathcal{T}=\bigcup_{i=1}^k\{A: \hbox{ $A$ is a $\mathcal{G}$-factor of $S_i$}\}.$$
 (We say that $\mathcal{T}$ is \emph{induced} by the flag $S_1\subseteq S_2\subseteq\dots\subseteq S_k$.)
\item Let $A_1,\dots,A_k$ be the minimal elements of $\mathcal{T}$, then they are all the $\mathcal{G}$-factors of
a certain element in $\mathcal{S}$, and for each $1\le i\le k$, the set $\{A\in \mathcal{T}: A\supsetneq A_i\}$ is also $\mathcal{G}$-nested defined
by induction.
\end{enumerate}
\end{df}

The wonderful compactification is defined as follows:

\begin{df}\label{df:wonderful compactification} Denote $Y^\circ=Y\setminus\cup_{G\in\mathcal{G}}G$. There is a natural locally closed embedding
$$Y^\circ\hookrightarrow Y\tm \prod_{G\in\mathcal{G}}Bl_GY.$$
The closure of this embedding, denoted by $Y_\mathcal{G}$, is called the \emph{wonderful compactification}
of $\mathcal{G}$.
\end{df}

The wonderful compactification $Y_\mathcal{G}$ of $\mathcal{G}$ has the following properties, where (1) and (2) are in Theorem 1.2 in \cite{Li} and (3) is clear from the proof there.
\begin{thm}\label{main wonder}
  The variety $Y_\mathcal{G}$ is nonsingular. For each $G\in\mathcal{G}$
  there is a nonsingular divisor $D_G$ on $Y_\mathcal{G}$ such that:
\begin{enumerate}
\item[(1)] The union of the divisors $D_G$ is  $Y_\mathcal{G}\setminus
Y^\circ$,
\item[(2)] Any collection of  the divisors $D_G$ intersects transversally. An
  intersection of divisors $D_{T_1}\cap\cdots\cap D_{T_r}$ is nonempty
  exactly when $\{T_1,\cdots,T_r\}$ forms a $\mathcal{G}$-nest.
\item[(3)] Each $D_G$ is the unique connected component of $\pi^{-1}(G)$ that maps surjectively to the subvariety $G$, where $\pi$ is
the natural morphism $Y_\mathcal{G}\to Y$. (This $D_G$ is called the \emph{dominant transform} of $G$ and denoted by $\widetilde{G}$ in \cite{Li}.)

\end{enumerate}
\end{thm}

The dominant transform can also be defined as follows. Let $\pi: \widetilde{Y}\to Y$ be the blow-up along a nonsingular subvariety $G\subsetneq Y$.
For any irreducible subvariety $V$ in $Y$, we define the dominant transform  of $V$, denoted by $\widetilde{V}$ or $V^\sim$, to be the strict
transform of $V$ when $V\nsubseteq G$, and $\pi^{-1}(V)$ when $V\subseteq G$. For a sequence of $N$ blow-ups $Y_N\to Y_{N-1}\to\cdots \to Y_1\to Y_0$
and a subvariety $V\subseteq Y_0$ we define the dominant transform $\widetilde{V}\subseteq Y_N$ (or denoted by $V^\sim$) to be the $N$-th iterated
dominant transform $(\cdots((V^\sim)^\sim)\cdots)^\sim$.

\medskip
It is known (cf. \cite{Li}) that $Y_{\mathcal{G}}$ can be constructed by a sequence of blow-ups as follows. Let $Y$ be a nonsingular variety,
$\mathcal{S}$ be an arrangement of subvarieties and $$\mathcal{G}=\{G_1,\dots, G_N\}$$ be a building set with respect to $\mathcal{S}$. Suppose the
subvarieties in $\mathcal{G}=\{G_1,\dots, G_N\}$ are indexed in an order compatible with inclusion relations, i.e. $i\le j $ if $ G_i\subseteq G_j$.
We define the triple $(Y_k, \mathcal{S}^{(k)}, \mathcal{G}^{(k)})$ inductively with respect to $k$, where $Y_k$ is a nonsingular variety,
$\mathcal{S}^{(k)}$ is an arrangement of subvarieties of $Y_k$ and $\mathcal{G}^{(k)}$ is a building set with respect to $\mathcal{S}^{(k)}$:

(1) For $k=0$, define $Y_0=Y$, $\mathcal{S}^{(0)}=\mathcal{S}$, $\mathcal{G}^{(0)}=\mathcal{G}=\{G_1,\dots, G_N\}$, $G^{(0)}_i=G_i$ for $1\le i\le
N$.

(2) Assume the triple $(Y_k, \mathcal{S}^{(k)}, \mathcal{G}^{(k)})$ has been constructed. Define $Y_{k}$ to be the blow-up of $Y_{k-1}$ along the nonsingular subvariety $G^{(k-1)}_k$. Define $G^{(k)}$ to be the dominant transform $(G^{(k-1)})^\sim$ for all $G\in
\mathcal{G}$.  Then  $\mathcal{G}^{(k)}:=\{G^{(k)}\}_{G\in\mathcal{G}}$ is a building set (by \cite{Li} Proposition 2.8). We denote the induced arrangement by $\mathcal{S}^{(k)}$.

(3) Continue the inductive construction until $k=N$. We get a nonsingular variety $Y_N$ and all elements in the building set $\mathcal{G}^{(N)}$ are
divisors. The resulting variety is isomorphic to $Y_\mathcal{G}$.

\medskip
For any $\mathcal{G}$-nest $\mathcal{T}$, define $$Y_k\mathcal{T}=\bigcap_{G\in \mathcal{T}} G^{(k)}.$$
The following property of $Y_k\mathcal{T}$ is
used often throughout the paper.
\begin{prop}\label{how strata change}
    Let $0\le k\le N-2$ and let $\mathcal{T}\subseteq\{G_{k+2},G_{k+3},\dots, G_N\}$ be a
    $\mathcal{G}$-nest. Then $Y_{k+1}\mathcal{T}$ is an irreducible nonsingular subvariety of
    $Y_{k+1}$ with the following property:

If $\{G_{k+1}\}\cup\mathcal{T}$ is not a $\mathcal{G}$-nest, then $G^{(k)}_{k+1}\cap Y_k\mathcal{T}=\emptyset$ and $Y_{k+1}\mathcal{T}\cong
Y_k\mathcal{T}$; otherwise, the intersection $G_{k+1}^{(k)}\cap Y_k\mathcal{T}$ is clean, $Y_{k+1}\mathcal{T}$ is isomorphic to the blow-up of
$Y_k\mathcal{T}$ along $G^{(k)}_{k+1}\cap Y_k\mathcal{T}$ with exceptional divisor $G^{(k+1)}_{k+1}\cap Y_{k+1}\mathcal{T}$ (where the intersection
is transverse), and the codimension of $G^{(k)}_{k+1}\cap Y_k\mathcal{T}$ in $Y_k\mathcal{T}$ equals to
$$\left\{
           \begin{array}{ll}
            \dim\cap_{G_{k+1}\subsetneq G\in\mathcal{T}}G-\dim {G_{k+1}}, & \hbox{ if $\{G:G_{k+1}\subsetneq G\in \mathcal{T}\}
                \neq \emptyset$;} \\
             \dim Y-\dim G_{k+1}, & \hbox{ otherwise.}
           \end{array}
         \right.$$

\end{prop}

\begin{proof}
We prove the statement by induction on $k$. The case $k=0$ is obvious. Now assume that the statement is true for $k$.

(i) Suppose that $\{G_{k+1}\}\cup\mathcal{T}$ is not a $\mathcal{G}$-nest. We will show that $G^{(k)}_{k+1}\cap Y_k\mathcal{T}=\emptyset$. As a
consequence we have $Y_{k+1}\mathcal{T}\cong Y_k\mathcal{T}$, since the center of the blow-up is away from $Y_k\mathcal{T}$.

We prove by contradiction. Assume  $G^{(k)}_{k+1}\cap Y_k\mathcal{T}\neq\emptyset$. Since $\mathcal{T}$ is a $\mathcal{G}$-nest,
$\{G^{(k)}\}_{G\in\mathcal{T}}$ is a $\mathcal{G}^{(k)}$-nest by \cite{Li} Proposition 2.8 (3). By Definition \ref{nest} (1), the nest
$\{G^{(k)}\}_{G\in\mathcal{T}}$ is induced by a flag $$S'_1\subseteq S'_2\subseteq\cdots\subseteq S'_l$$ where $S'_i\in \mathcal{S}^{(k)}$. We claim
that $\{G^{(k)}_{k+1}\} \cup \{G^{(k)}\}_{G\in \mathcal{T}}\subseteq \mathcal{G}^{(k)}$ is a $\mathcal{G}^{(k)}$-nest induced by the flag
\begin{equation}\label{e-flag}\big{(}G^{(k)}_{k+1}\cap S'_1\big{)}\subseteq S'_1\subseteq S'_2\subseteq\cdots\subseteq
S'_l.
\end{equation}
Indeed, since $Y_k\mathcal{T}=S'_1$, we know $G^{(k)}_{k+1}\cap S'_1\neq\emptyset$. By \cite{Li} Lemma 2.4 (ii), the $\mathcal{G}^{(k)}$-factors of
$G^{(k)}_{k+1}\cap S'_1$ are $G^{(k)}$ and some $\mathcal{G}^{(k)}$-factors of $S'_1$, hence our claim follows. Then \cite{Li} Proposition 2.8 (3)
asserts that, since $\{G^{(k)}_{k+1}\} \cup \{G^{(k)}\}_{G\in \mathcal{T}}$ is a $\mathcal{G}^{(k)}$-nest, $\{G_{k+1}\}\cup\mathcal{T}$ must be a
$\mathcal{G}$-nest. But by assumption $\{G_{k+1}\}\cup\mathcal{T}$ is not a $\mathcal{G}$-nest, contradition.

(ii) Suppose that $\mathcal{T}\cup \{G_{k+1}\}$ is a
$\mathcal{G}$-nest. Let the $\mathcal{G}^{(k)}$-factors of $Y_k\mathcal{T}$ be
$G'_1,\dots, G'_r$. Then they are minimal elements in the
$\mathcal{G}^{(k)}$-nest $\{G^{(k)}\}_{G\in\mathcal{T}}$, by the
definition of nest.
Assume without loss of generality that the first $m$ subvarieties $G'_1,\dots,G'_m$
contain $G^{(k)}_{k+1}$. Define $A=\cap_{i=1}^mG'_i$,
$B=\cap_{i=m+1}^rG'_i$, then $Y_k\mathcal{T}=A\cap B$ is the
$G^{(k)}_{k+1}$-factorization of $Y_k\mathcal{T}$ by \cite{Li} Definition-Lemma 2.6.

Notice that for $p,q\ge k+2$ and $G_p^{(k)}\subseteq G_q^{(k)}$, we
have $G_p^{(k+1)}\subseteq G_q^{(k+1)}$ because strict transforms
keep the inclusion relation. Moreover, since $G'_1,\dots, G'_r$ are the
minimal elements in $\mathcal{G}^{(k)}$ which contain
$Y_k\mathcal{T}$,
the subvariety $Y_{k+1}\mathcal{T}$ is the intersection $\cap_{i=1}^r\widetilde{G}'_i$. Then
$$\widetilde{A}=\bigcap_{i=1}^m\widetilde{G}'_i,\quad
\widetilde{B}=\bigcap_{i=m+1}^r\widetilde{G}'_i,\quad
(A\cap B)^\sim =\widetilde{A}\cap \widetilde{B}=\bigcap_{i=1}^m\widetilde{G}'_i$$
by \cite{Li} Lemma 2.9. Thus $Y_{k+1}\mathcal{T}=({Y_{k}\mathcal{T}})^\sim$. By the definition of arrangement we know that
$Y_k\mathcal{T}$ and $G^{(k)}_{k+1}$ intersect cleanly, so
$Y_{k+1}\mathcal{T}$ is the blow-up of $Y_k\mathcal{T}$ along the
center $G^{(k)}_{k+1}\cap Y_k\mathcal{T}$. The exceptional divisor
is the preimage of the center, hence is $
G^{(k+1)}_{k+1}\cap Y_{k+1}\mathcal{T}$. Since $G^{(k+1)}_{k+1}$ and $Y_{k+1}\mathcal{T}$ intersect cleanly and since the divisor $G^{(k+1)}_{k+1}$ does not contain $Y_{k+1}\mathcal{T}$, we can see that the intersection  $G^{(k)}_{k+1}\cap Y_k\mathcal{T}$ is actually transversal.

The codimension of the center $Y_k\mathcal{T}\cap G^{(k)}_{k+1}$ in
$Y_k\mathcal{T}$ equals
$${\rm codim}({A\cap B\cap G^{(k)}_{k+1}}, A\cap B)={\rm codim}({G^{(k)}_{k+1}\cap B}, A\cap B)={\rm codim}({G^{(k)}_{k+1}}, A),$$
where the second equality is because of the transversality of the intersection $G^{(k)}_{k+1}\cap B$. If no elements in $\mathcal{T}$ contain
$G_{k+1}$, then $A=Y$ and $${\rm codim}({G^{(k)}_{k+1}}, A)=\dim Y-\dim G_{k+1};$$ otherwise
$$\begin{array}{ll}
{\rm codim}({G^{(k)}_{k+1}}, A)&={\rm codim}({G^{(k)}_{k+1}},\bigcap_{i=1}^mG'_i)={\rm codim}({G_{k+1}}, \bigcap_{G_{k+1}\subsetneq
G\in\mathcal{T}}G)\\&=\dim\bigcap_{G_{k+1}\subsetneq G\in\mathcal{T}}G-\dim {G_{k+1}}. \end{array}$$
 Thus the proof is complete.
\end{proof}

%%%%%%%%%%%%%%%%%%%%%%%%%%%%%%%%%%%%%%%%%%%%%%%%%%%%%%%%%%%%%%%%%%%%%%%555
\ss{The motive of wonderful compactifications}\label{main thm setting}

\noindent\textbf{Notations}:

\smallskip

\noindent $\bullet$\, Let $Y$ be a nonsingular quasi-projective variety with an arrangement of subvarieties $\mathcal{S}$. Let $\mathcal{G}$ be a
building set with respect to $\mathcal{S}$. Let $Y_\mathcal{G}$ be the wonderful compactification. Let $\mathcal{T}$ be a $\mathcal{G}$-nest.

\smallskip
\noindent $\bullet$\, For $T\in\mathcal{G}$, define $D_T$ to be the divisor $T^{(N)}$ in $Y_\mathcal{G}$. When no confusion arise, we use the same
notation $D_T$ for its restriction to a subvariety of $Y_\mathcal{G}$.

\smallskip
\noindent $\bullet$\, Denote $j_{\mathcal{T}}: Y_\mathcal{G}\mathcal{T}\to Y_\mathcal{G}$ to be the natural imbedding. Denote $g_{\mathcal{T}}:
Y_\mathcal{G}\mathcal{T}\to Y_0\mathcal{T}$ to be the restriction of the natural morphism $Y_\mathcal{G}\to Y$.

\smallskip
\noindent $\bullet$\, Suppose $j: B\to C$ and $g: B\to D$ are two morphisms of varieties. Denote by $(j\btm g): B\to C\tm D$  the composition of the
diagonal map $\Delta$ with $f\tm g$:
$$(j\btm g): B\stackrel{\Delta}{\to} B\tm B\stackrel{f\tm g}{\to} C\tm D.$$

\smallskip
\noindent $\bullet$\, Given $a\in A(P)$, denote by $\{a\}_i$ the image of the projection $A(P)\to A^i(P)$ of the Chow ring to its degree $i$ direct
summand, i.e., taking the codimension $i$ part of $a$.

\smallskip
\noindent $\bullet$\, We set $\bigcap_{G\subsetneq T\in\mathcal{T}}T=Y$ if no $T$ satisfies $G\subsetneq T\in\mathcal{T}$. Define
$$r_G:=\dim(\bigcap_{G\subsetneq T\in\mathcal{T}}T)-\dim G.$$
Define $$N_G:=N_G(\bigcap_{G\subsetneq T\in\mathcal{T}}T)|_{Y_0\mathcal{T}},$$ the restriction to $Y_0\mathcal{T}$ of the normal bundle of $G$ in the
ambient space $(\bigcap_{G\subsetneq T\in\mathcal{T}}T)$. Define $$M_\mathcal{T}:=\big{\{}\underline{\mu}=\{\mu_G\}_{G\in\mathcal{G}}: 1\le \mu_G\le
r_G-1, \mu_G\in\mathbb{Z}\big{\}}$$ and define $||\underline{\mu}||:=\sum_{G\in\mathcal{G}}\mu_G$ for $\underline{\mu}\in M_\mathcal{T}$.

\begin{thm}\label{main thm wonderful} We have the Chow group decomposition
$$A^*Y_\mathcal{G}=A^*Y\oplus\bigoplus_{\mathcal{T}}\bigoplus_{\underline{\mu}\in
M_\mathcal{T}} A^{*-\|\underline{\mu}\|}(Y_0\mathcal{T})$$ where $\mathcal{T}$ runs through all
$\mathcal{G}$-nests.

Moreover, when $Y$ is complete, we have the Chow motive decomposition
$$h(Y_\mathcal{G})=h(Y)\oplus\bigoplus_{\mathcal{T}}\bigoplus_{\underline{\mu}\in
M_\mathcal{T}} h(Y_0\mathcal{T})(||\underline{\mu}||)$$ where $\mathcal{T}$ runs through all
$\mathcal{G}$-nests.
\end{thm}

\begin{thm}\label{main thm correspondence}
The correspondence that gives each of the above direct summand can be explicitly expressed as follows,
\begin{align*}
\alpha&: h(Y_\mathcal{G})\to h(Y_0\mathcal{T})(\|\underline{\mu}\|)\\
 \alpha&=(j_{\mathcal{T}}\btm g_{\mathcal{T}})_*\prod_{G\in \mathcal{T}} \bigg{\{}
c\Big{(}g_{\mathcal{T}}^*(N_G)\ts \mathcal{O}\big{(}-\sum_{(\bigstar)}D_{G'}\big{)}\Big{)}
\frac{1}{1+D_G}\bigg{\}}_{r_G-1-\mu_G},
\end{align*}
here $c$ is total Chern class, the subscript ${r_G-1-\mu_G}$ means the codimension ${r_G-1-\mu_G}$ part, and condition $(\bigstar)$ is: $
G'\subsetneq G$ and $\mathcal{T}\cup\{G'\}$ is a $\mathcal{G}$-nest.

The inverse correspondence is
\begin{align*}\beta&: h(Y_0\mathcal{T})(\|\underline{\mu}\|) \to h(Y_\mathcal{G})\\
\beta&=(g_{\mathcal{T}}\btm j_{\mathcal{T}})_*\prod_{G\in \mathcal{T}}\big{(}-D_G\big{)}^{\mu_G-1}.
\end{align*}

\end{thm}

%%%%%%%%%%%%%%%%%%%%%%%%%%%%%%%%%%%%%%%%%%%%%%%%%%%%%%%%%%%%%%%%%%%%%%%

\sss{Proof of the Theorem \ref{main thm wonderful}}

\begin{lem}\label{one blow-up}
Given a $\mathcal{G}$-nest $\mathcal{T}\subseteq\{G_{k+2},\dots,G_N\}$. Suppose
$\mathcal{T}':=\mathcal{T}\cup\{G_{k+1}\}$ is also a $\mathcal{G}$-nest. Define $r=r_{k,\mathcal{T}}$ to be $$\left\{
           \begin{array}{ll}
            \dim\cap_{G_{k+1}\subsetneq G\in\mathcal{T}}G-\dim {G_{k+1}}, & \hbox{ if $\{G:G_{k+1}\subsetneq G\in \mathcal{T}\}
                \neq \emptyset$;} \\
             \dim Y-\dim G_{k+1}, & \hbox{ otherwise.}
           \end{array}
         \right.$$
Then the following Chow group decomposition holds:
$$A^*(Y_{k+1}\mathcal{T})=A^*(Y_{k}\mathcal{T})\oplus
\bigoplus_{t=1}^{r-1}A^{*-t}(Y_k\mathcal{T}').$$ When $Y$ is complete, we also have the motivic
decomposition
$$h(Y_{k+1}\mathcal{T})=h(Y_{k}\mathcal{T})\oplus
\bigoplus_{t=1}^{r-1} h(Y_k\mathcal{T}')(t).$$
\end{lem}

\begin{proof}
Apply the well known blow-up formula for the Chow group and for the Chow motive (Theorem \ref{motive for a blow-up}) to Proposition
 \ref{how strata change} immediately gives the conclusion.
\end{proof}

Iteratively applying the above lemma gives the proof of Theorem \ref{main thm
wonderful}:

\begin{proof}[Proof of Theorem \ref{main thm wonderful}]

Define $$M^{(k)}_\mathcal{T}=\big{\{}\underline{\mu}=\{\mu_G\}_{G\in\mathcal{G}}: 1\le \mu_G\le
\dim (\bigcap_TT^{(k)})-\dim
G^{(k)}-1, \mu_G\in\mathbb{Z}\big{\}}$$ where $T$ runs through the subvarieties in $
\mathcal{T}$ such that $G^{(k)}\subsetneq T^{(k)}$. Define $||\underline{\mu}||:=\sum_{G\in\mathcal{G}}\mu_G$ for
$\underline{\mu}\in M^{(k)}_\mathcal{T}$.

We prove the following statement using a downward induction on $k$ :
\emph{
\begin{equation}\label{equation}
A^*Y_\mathcal{G}=A^*Y_k\oplus\bigoplus_{\mathcal{T}}\bigoplus_{\underline{\mu}\in
M^{(k)}_\mathcal{T}} A^{*-\|\underline{\mu}\|}(Y_k\mathcal{T}).
\end{equation}
where $\mathcal{T}$ runs through all $\mathcal{G}$-nest such that
$\mathcal{T}\subseteq\{G_{k+1},G_{k+2},\dots, G_N\}$.}
\medskip

The assertion for $k=N$ is trivial because all $G^{(N)}$ are divisors in $Y_\mathcal{G}$ hence of
codimension 1 and $M^{(k)}_\mathcal{T}=\emptyset$.

Assume (\ref{equation}) has been proved for $k+1$, i.e.,
$$A^*Y_\mathcal{G}=A^*Y_{k+1}\oplus\bigoplus_{\mathcal{T}}\bigoplus_{\underline{\mu}\in
M^{(k+1)}_\mathcal{T}} A^{*-\|\mu\|}(Y_{k+1}\mathcal{T})$$ where $\mathcal{T}$ runs through all
$\mathcal{G}$-nest such that $\mathcal{T}\subseteq\{G_{k+2},G_{k+3},\dots, G_N\}$. Apply Lemma
\ref{one blow-up}, we have
\begin{align}\label{eq:chow}
A^*Y_{\mathcal{G}}=&A^*Y_k\bigoplus
\bigg{(}\bigoplus_{t=1}^{\textrm{codim}(G_{k+1},Y)-1}A^{*-t}(G^{(k)}_{k+1})\bigg{)}\\
&\bigoplus\bigg{(}\bigoplus_{\mathcal{T}}\bigoplus_{\underline{\mu}\in M^{(k+1)}_{\mathcal{T}}} A^{*-\|\mu\|}(Y_k\mathcal{T})
\bigg{)}\nonumber\\
&\bigoplus\bigg{(}\bigoplus_{\mathcal{T}}\bigoplus_{\underline{\mu}\in
M^{(k+1)}_{\mathcal{T}}}\bigoplus_{t=1}^{r_{k+1,\mathcal{T}}-1}A^{*-\|\mu\|-t}\big{(}Y_k(\{G_{k+1}\}\cup\mathcal{T})\big{)}\bigg{)}.\nonumber
\end{align}
 This immediately gives the Chow group decomposition (\ref{equation}) for
$k$. Indeed, any $\mathcal{G}$-nest contained in $\{G_{k+1},G_{k+2},\dots, G_N\}$ must be one of the three: $\{G_{k+1}\}$, a $\mathcal{G}$-nest
$\mathcal{T}$ contained in $\{G_{k+2},G_{k+3},\dots, G_N\}$, or $\{G_{k+1}\}\cup\mathcal{T}$. They correspond to the second, third and last summands
in (\ref{eq:chow}) respectively. (Notice that $Y_k(\{G_{k+1}\}\cup \mathcal{T}))=\emptyset$ if $\{G_{k+1}\}\cup \mathcal{T}$ is not a
$\mathcal{G}$-nest by Proposition \ref{how strata change}.)

Therefore, the Chow group decomposition (\ref{equation}) holds for all $k$, in particular the case
$k=0$ gives the desired Chow group decomposition. For the proof of the Chow motive decomposition, we can either repeat the above proof almost word by word or, as the referee pointed out, notice that the Chow motive decomposition follows from the result on the Chow groups and Manin's identity principle.
\end{proof}

\sss{Proof of Theorem \ref{main thm correspondence}}

First, we introduce some notations. For a given $\mathcal{G}$-nest $\mathcal{T}$,

\noindent$\bullet $\, Define $\mathcal{T}_k:=\mathcal{T}\bigcap\{G_{k+1},G_{k+2},\dots,G_N\}$ for $0\le k\le N$. Then we have a chain of
$\mathcal{G}$-nests $\mathcal{T}_0\supseteq\mathcal{T}_1\supseteq\dots\supseteq\mathcal{T}_N$, where $\mathcal{T}_0=\mathcal{T}$ and
$\mathcal{T}_N=\emptyset$.

\noindent$\bullet $\,  For $\underline{\mu}\in \mathcal{M}_\mathcal{T}$ and $1\le i\le N$, define $$\mu_i:=\left\{
                                        \begin{array}{ll}
                                          \mu_{G_i}, & \hbox{if $G_i\in\mathcal{T}$;} \\
                                          0, & \hbox{otherwise.}
                                        \end{array}
                                      \right.$$

\noindent$\bullet $\,  $j_{kl}$ and $g_{kl}$ ($N\ge k>l\ge 0$) are the natural morphisms as in the following diagram
\begin{equation}\xymatrix{
 Y_N\mathcal{T}_0\quad \ar`u[r]`[rrrr]^{j_{\mathcal{T}}}[rrrr] \ar`l[d]`[dddd]_{g_{\mathcal{T}}}[dddd]
\ar[r]^{j_{N0}}\ar[d]_{g_{N0}}\ar@{}[dr]|{} & Y_N\mathcal{T}_1 \ar[r]^-{j_{N1}}\ar[d]_{g_{N1}}\ar@{}[dr]|{}& ...\ar[r]
&Y_N\mathcal{T}_{N-1}\ar[r]^-{j_{N,N-1}}\ar[d]_{g_{N,N-1}} &
 Y_N\mathcal{T}_N\ar@{.>}[ld]_{\alpha_{N}}\\
Y_{N-1}\mathcal{T}_0 \ar[r]^{j_{N-1,0}}\ar[d]_-{g_{N-1,0}} \ar@{}[dr]|{}& Y_{N-1}\mathcal{T}_1
\ar[r]^-{j_{N-1,1}}\ar[d]_-{g_{N-1,1}}& ...\ar[r] & Y_{N-1}\mathcal{T}_{N-1}\ar@/_/@{.>}[ru]_{\beta_{N}}\\
...\ar[d]_-{g_{20}}\ar@{}[dr]|{}& ...\ar[d]_-{g_{21}}& ...\\
Y_1\mathcal{T}_0\ar[r]^{j_{10}}\ar[d]_{g_{10}}&Y_1\mathcal{T}_1\ar@{.>}[ld]_{\alpha_1}\\
Y_0\mathcal{T}_0\ar@/_/@{.>}[ru]_{\beta_1}} \label{big_diagram}
\end{equation}

\begin{lem}\label{g^-1 divisor}
Denote by $g: Y_k\to Y_{k-1}$ the natural morphism. Then for $l\le k-1$, we have
$$g^{-1}(G^{(k-1)}_l)=G^{(k)}_l.$$
\end{lem}
\begin{proof}
First, we claim that $G^{(k-1)}_l\nsupseteq G^{(k-1)}_k$. Otherwise $G_l\supseteq G_k$ since they are the respective images of $G^{(k-1)}_l$ and $G^{(k-1)}_k$ under
$Y_k\to Y_0$. But then by the assumption that the order of $\{G_i\}$ is compatible with inclusion relations, we obtain a contradiction $l\ge k$.

Next, it is easy to see that $G^{(k-1)}_l\nsubseteq G^{(k-1)}_k$ since $G^{(k-1)}_l$ is a divisor. Now we know that the two nonsingular subvarieties $G^{(k-1)}_l$ and $G^{(k-1)}_k$ intersect cleanly and neither one contains the other, therefore they
must intersect transversally. Then it is standard to show by local coordinates calculation that the
following isomorphism between ideal sheaves holds:
$$g^{-1}\mathcal{I}(G^{(k-1)}_l)\cdot\O_{Y_k}\cong \mathcal{I}(G^{(k)}_l).$$ The desired conclusion follows from this.
\end{proof}

\begin{lem}\label{property of big diagram}
In Diagram (\ref{big_diagram}), all squares are fiber squares. Moreover, for any $N\ge k>l\ge 0$, we have
\begin{enumerate}
\item[(i)] $j_{kl}$ is injective;
\item[(ii)] If $G_k\in\mathcal{T}$, then $g_{kl}$ is the projection of a projective bundle with fiber isomorphic to a projective space of dimension
${r_{k,\mathcal{T}}-1}$;
\item[(iii)] If $G_k\notin\mathcal{T}$ but $\{G_k\}\cup\mathcal{T}_l$ is a
$\mathcal{G}$-nest, then $g_{kl}$ is the blow-up of $Y_{k-1}\mathcal{T}_l$ along the center $G^{(k-1)}_k\cap
Y_{k-1}\mathcal{T}_l$;
\item[(iv)] If $\{G_k\}\cup\mathcal{T}_l$ is not a $\mathcal{G}$-nest, then $g_{kl}$ is an isomorphism.
\end{enumerate}
\end{lem}
\begin{proof}
It is obvious that $j_{kl}$ is injective.

Now we show that $g_{kl}$ is the projection of a projective bundle if $G_k\in\mathcal{T}$. By
Proposition \ref{how strata change}, the variety $Y_k\mathcal{T}_k$ is the blow-up of $Y_{k-1}\mathcal{T}_k$
along the center $Y_{k-1}\mathcal{T}_{k-1}$, and the exceptional divisor is $Y_k\mathcal{T}_{k-1}$ (note that $Y_{k-1}\mathcal{T}_k\cap G^{(k-1)}_k=Y_{k-1}\mathcal{T}_{k-1}$ and $Y_k\mathcal{T}_k\cap G^{(k)}_k=Y_k\mathcal{T}_{k-1}$).
Therefore $g_{k,k-1}: Y_k\mathcal{T}_{k-1}\to Y_{k-1}\mathcal{T}_{k-1}$ is a projective bundle, and
the dimension of a fibre is $r_{k,\mathcal{T}}-1$. Next we show that for any $l\le k-1$, $g_{kl}$
is the restriction of $g_{k,k-1}$ to a smaller base $Y_{k-1}\mathcal{T}_l$, which will then show that $g_{kl}$ is also a
projective bundle with fiber of the same dimension $r_{k,\mathcal{T}}-1$. Fix $k$ and use downward
induction on $l$. By inductive assumption, $g_{k,l+1}$ is a restriction of $g_{k,k-1}$. Since
$$g_{k,l+1}^{-1}(G^{(k-1)}_{l+1}\cap Y_{k-1}\mathcal{T}_l)=G^{(k)}_{l+1}\cap Y_{k}\mathcal{T}_l$$
by Lemma \ref{g^-1 divisor}, the restriction of the projective bundle $g_{k,l+1}$ to a smaller base
space $Y_{k-1}\mathcal{T}_{l}=Y_{k-1}\mathcal{T}_{l+1}\cap G^{(k-1)}_{l+1}$ is exactly $g_{kl}$.

Next, we show $g_{kl}$ is birational if $G_k\notin\mathcal{T}$. This is again implied by
Proposition \ref{how strata change}. Notice that $G^{(k-1)}_k$ is minimal in
$$\mathcal{T}':=\{G^{(k-1)}_k\}\cup \{G^{(k-1)}\}_{G\in\mathcal{T}_l}.$$
If $\mathcal{T}'$ is a $\mathcal{G}^{(k-1)}$-nest, then $g_{kl}:Y_{k}\mathcal{T}_l\to
Y_{k-1}\mathcal{T}_l$ is a blow-up along the center $G^{(k-1)}_k\cap Y_{k-1}\mathcal{T}_l$;
otherwise, $g_{kl}$ is an isomorphism. In both cases, $g_{kl}$ is birational.

Finally, all squares in Diagram (\ref{big_diagram}) are fiber squares since $\forall l\le k-2$,
$g_{kl}$ is a restriction of $g_{k,l+1}$. The proof is complete.
\end{proof}

The following lemma computes the composition of correspondences in certain diagrams. The author thanks the referee to suggest a proof much simpler
than the original proof given by the author.
\begin{lem}\label{comp} Let $W,U,V,X,Y,Z$ be nonsingular quasi-projective varieties. Suppose the square in the following diagram is a fiber square.
$$ \xymatrix{
 W \ar[r]^{j_{3}}\ar[d]_{g_{3}} \ar@{}[dr]|{\square}& U\ar[d]_{g_{2}}\ar[r]^{j_{2}} & X \ar@/_/@{.>}[ld]_{\alpha_{2}}\\
V \ar[r]^{j_{1}}\ar[d]_{g_{1}} & Y\ar@/_/@{.>}[ur]_{\beta_2}\ar@/_/@{.>}[ld]_{\alpha_1}\\
Z\ar@/_/@{.>}[ru]_{\beta_1}}$$ and suppose that $\dim W-\dim V=\dim U-\dim Y$ and that $j_k, g_k (1\le k\le 3)$ are proper. Take $\gamma_1, \gamma_1'\in A(V),
\gamma_2, \gamma_2'\in A(U)$ and define correspondences
$$\alpha_k=(j_{k},g_{k})_*\gamma_k, \quad \beta_k=(g_{k},j_{k})_*\gamma'_k, \quad \textrm{for $k=1,2$.}$$
Then we have
\begin{align}
\label{eq:alpha}\alpha_1\alpha_2&=(j_2j_3,g_1g_3)_*(j_{3}^*\gamma_2\cdot g_{3}^*\gamma_1),\\
\label{eq:beta}\beta_2\beta_1&=(g_1g_3,j_2j_3)_*(g_{3}^*\gamma'_1\cdot j_{3}^*\gamma'_2).
\end{align}
\end{lem}

\begin{proof}
By abuse of notation, for $\gamma\in A(V)$ we use the same $\gamma$ to denote the correspondence $(\Delta_V)_*(\gamma)\in A(V\times V)$ where $\Delta_V: V\to V\times V$ is the diagonal embedding. For a map $j:U\to X$, we denote by $j_*$ the correspondence $\Gamma_j$ (i.e. the graph of $j$) and by $j^*$ the correspondence $\Gamma_j'$ (i.e. the transpose of $\Gamma_j$).

First observe that $\alpha_k=g_{k*}\circ\gamma\circ j_k^*$ for $k=1,2.$ Indeed, by properties of correspondences (cf. \cite{Fulton} Prop 16.1.1(c)),
we have $\Gamma_j\circ \gamma=(1_U\times j)_*\gamma$, $\gamma\circ\Gamma_g'=(g\times 1_U)_*\gamma,$ so $$g_{k*}\circ\gamma_k\circ
j_k^*=\Gamma_{g_k}\circ\gamma_k\circ\Gamma_{j_k}^*=(g_k\times j_k)_*\gamma_k=(g_k,j_k)_*\gamma_k=\alpha_k\quad \hbox{for } k=1,2.$$ With the above
observation, (\ref{eq:alpha}) is equivalent to
$$g_{1*}\gamma_1j^*_1g_{2*}\gamma_2j^*_2=(g_1g_3)_*(j^*_3\gamma_2\cdot g^*_3\gamma_1)(j_2j_3)^*.$$
So it suffices to prove
\begin{equation}\label{eq:correspondences}\gamma_1j^*_1g_{2*}\gamma_2=g_{3*}(j^*_3\gamma_2\cdot g^*_3\gamma_1)j_3^*.
\end{equation}
For any $u\in A(U)$, we have
$$\gamma_1j^*_1g_{2*}\gamma_2(u)=\gamma_1g_{3*}j_3^*\gamma_2(u)=g_{3*}(g^*_3\gamma_1\cdot j^*_3(\gamma_2u))=g_{3*}(g^*_3\gamma_1\cdot j^*_3\gamma_2) j^*_3(u)$$
where the first ``='' is because of $\dim W-\dim V=\dim U-\dim Y$, the second ``='' is because of the projection formula. Then we apply Manin's
Identity Principle to obtain (\ref{eq:correspondences}), hence (\ref{eq:alpha}). The identity (\ref{eq:beta}) can be obtained by transposing
(\ref{eq:alpha}).
\end{proof}

Now we state a simple lemma and omit the proof.

\begin{lem}\label{lm:simple} If $A$, $B_i$, $C_{ij}$ are motives such that
\begin{enumerate}
\item[(i)] $\bigoplus_i\alpha_i:\, A\cong\bigoplus_i B_i$ is an isomorphism with inverse $\sum_i\beta_i$, and
\item[(ii)] $\bigoplus_j\alpha_{ij}:\, B_i\cong\bigoplus_j C_{ij}$ is an isomorphism with inverse
$\sum_j\beta_{ij}$,
\end{enumerate}
then the correspondence $\bigoplus_{i,j}\alpha_{ij}\circ\alpha_i$ gives an isomorphism  $A\cong\bigoplus_{i,j}C_{ij}$ with inverse
$\sum_{i,j}\beta_i\circ\beta_{ij}$.
\end{lem}

For $G_k\in \mathcal{T}$, define $h_k\in A^1(Y_k\mathcal{T}_{k-1})$ to be first Chern class of the invertible sheaf $\mathcal{O}(1)$ of the
projective bundle $g_{k,k-1}$. Define
$$\alpha_k=\left\{
             \begin{array}{ll}
               (j_{k,k-1}\btm g_{k,k-1})_*1, & \hbox{if $G_k\notin\mathcal{T}$;} \\
               (j_{k,k-1}\btm g_{k,k-1})_*
              \Big{(}\big{\{} g_{k,k-1}^*c(N_k)\dfrac{1}{1-h_k}\big{\}}_{r_k-1-\mu_k}\Big{)},
               & \hbox{if $G_k\in\mathcal{T}$,}
             \end{array}
           \right.
$$
where $N_k:=N_{Y_{k-1}\mathcal{T}_{k-1}}Y_{k-1}\mathcal{T}_{k}$. Define
$$\beta_k=\left\{
             \begin{array}{ll}
               (g_{k,k-1}\btm j_{k,k-1})_*1, &  \hbox{if $G_k\notin\mathcal{T}$;} \\
               (g_{k,k-1}\btm j_{k,k-1})_*h_k^{\mu_k-1},
               & \hbox{if $G_k\in\mathcal{T}$.}
             \end{array}
           \right.
$$ Thanks to the blow-up formula
of motives (Theorem \ref{motive for a blow-up}), the correspondence
$$a_k: h(Y_k\mathcal{T}_k)(\sum_{i=k+1}^N\mu_i)\to
h(Y_{k-1}\mathcal{T}_{k-1})(\sum_{i=k}^N\mu_i)$$ expresses
$h(Y_{k-1}\mathcal{T}_{k-1})(\sum_{k}^N\mu_i)$ as a direct summand of
$h(Y_k\mathcal{T}_k)(\sum_{k+1}^N\mu_i)$ with right inverse $\beta_k$.

By Lemma \ref{lm:simple}, the correspondence  $$\alpha_{\mathcal{T},\underline{\mu}}: h(Y_\mathcal{G})\to
 h(Y_0\mathcal{T})(\|\underline{\mu}\|)$$
that gives the direct summand $h(Y_0\mathcal{T})(\|\underline{\mu}\|)$ in Theorem \ref{main thm
wonderful} can be expressed as the composition $\alpha_1\circ\alpha_2\circ\dots\circ\alpha_N$, with
right inverse $\beta_N\circ\dots\circ\beta_1$. Therefore we have

\begin{prop}\label{composition}
Denote by $f_k: Y_N\mathcal{T}_0\to Y_k\mathcal{T}_{k-1}$ the natural map in Diagram
(\ref{big_diagram}). (i.e. $g_{k+1,k-1}\dots\circ g_{N,k-1}\circ j_{N,k-2}\circ\dots\circ j_{N0}.$)
Then
\begin{align*}
\alpha_1\circ\dots\circ\alpha_N&=(j_{\mathcal{T}}\btm g_{\mathcal{T}})_*\prod_{G_k\in \mathcal{T}}
\big{\{}f_k^*g_{k,k-1}^*c(N_k)\frac{1}{1-f_k^*h_k}\big{\}}_{r_k-1-\mu_k}, \\
\beta_N\circ\dots\circ\beta_1&=(g_{\mathcal{T}}\btm j_{\mathcal{T}})_*\prod_{G_k\in
\mathcal{T}}f_k^*h_k^{\mu_k-1}.
\end{align*}
\end{prop}
\begin{proof}
  Combine Lemma \ref{property of big diagram} and Lemma \ref{comp} with the above
 discussion.
\end{proof}

The following two standard facts about normal bundles of subvarieties are used in the proof of
Theorem \ref{main thm correspondence}.
\begin{fact}\label{normal bundle_transversal}
Let $Z$ be a nonsingular variety. Let $Y, W$ be nonsingular proper subvarieties of $Z$ and assume $Y$ intersects transversally with $W$. Let $\pi:
\widetilde{Z}\to Z$ be the blow-up of $Z$ along $W$ and let $\widetilde{Y}$ be the strict transform of $Y$. Then
$$N_{\widetilde{Y}}\widetilde{Z}\simeq \pi^*N_YZ.$$
\end{fact}

\begin{fact}\label{normal bundle_subset}
Let $W\subsetneq Y\subsetneq Z$ be nonsingular varieties and $\pi: \widetilde{Z}\to Z$ be the
blow-up of $Z$ along $W$. Denote by $\Y$ the strict transform of $Y$,  and denote by $E$ the
exceptional divisor on $\widetilde{Y}$. Then
$$N_{\widetilde{Y}}\widetilde{Z}\simeq \pi^*N_YZ\otimes \O(-E).$$
\end{fact}

\begin{proof}[Proof of the above two facts]
Prove by local coordinates. Or see \cite{Fulton}.
\end{proof}

\begin{proof}[Proof of Theorem \ref{main thm correspondence}] To conclude Theorem \ref{main thm correspondence}
from Proposition \ref{composition}, we prove in three steps.
\medskip

{Step 1}: Show $f_k^*h_k=-D_{G_k}|_{Y_N\mathcal{T}_0}$.

\medskip
Recall that for $G_k\in \mathcal{T}$, $h_k$ is first Chern class of the invertible sheaf $\mathcal{O}(1)$ of the projective bundle $g_{k,k-1}$.

Consider the following diagram (not necessary a fiber square) where $\pi$ and $j$ are the natural
morphisms:
$$\xymatrix{Y_N\mathcal{T}_0\ar[r]^-{j_{\mathcal{T}}}\ar[d]^{f_k} & Y_N\ar[d]^{\pi}\\ Y_k\mathcal{T}_{k-1}
\ar[r]^-{j} & Y_k \,.}$$ By Proposition \ref{how strata change}, $Y_k\mathcal{T}_{k-1}$ is the exceptional divisor of the blow-up $g_{k,k-1}:
Y_k\mathcal{T}_{k-1}\to Y_{k-1}\mathcal{T}_{k-1}$, so $h_k=-j_{k,k-1}^*[Y_k\mathcal{T}_{k-1}]$. Since $Y_k\mathcal{T}_{k-1}$ is the transversal
intersection $Y_k\mathcal{T}_{k}\cap G^{(k)}_k$, $h_k=-j^*[G^{(k)}_k]$. Then
$$f_k^*h_k=-f_k^*j^*[G^{(k)}_k]=-j_{\mathcal{T}}^*\pi^*[G^{(k)}_k]=-j_{\mathcal{T}}^*D_{G_k}
=-D_{G_k}|_{Y_N\mathcal{T}_0}.$$ where the third equality can be proved by successively applying Lemma \ref{g^-1 divisor}.

\smallskip
{Step 2}: Let $0\le s<k\le N$. Denote $g_{sk}: Y_{s}\mathcal{T}_{k}\to Y_{s-1}\mathcal{T}_{k}$ to be the natural map induced from $Y_s\to Y_{s-1}$.
We claim the following:

If $G_k\in \mathcal{T}$ (hence $\mathcal{T}_{k-1}=\mathcal{T}_k\cup \{G_k\}$), then the normal bundle $N_{Y_{s}\mathcal{T}_{k-1}}Y_s\mathcal{T}_k$ is
isomorphic to $$\left\{
  \begin{array}{ll}
    g_{s,k-1}^*\big{(}N_{Y_{s-1}\mathcal{T}_{k-1}}Y_{s-1}\mathcal{T}_k\big{)}
\ts(-[G^{(s)}_s]|_{Y_s\mathcal{T}_{k-1}}), & \hbox{if (**) holds;} \\
\\
    g_{s,k-1}^*\big{(}N_{Y_{s-1}\mathcal{T}_{k-1}}Y_{s-1}\mathcal{T}_k\big{)}, & \hbox{otherwise.}
  \end{array}
\right. $$ where condition (**) is: $G_s\subsetneq G_k$ and  $\mathcal{T}_k\cup \{G_s\}$ is a
$\mathcal{G}$-nest.

\smallskip For the proof, we discuss three cases.

\smallskip
Case (i): condition (**) holds. It is a direct conclusion of Fact \ref{normal bundle_subset}. Indeed, to apply Fact \ref{normal bundle_subset} we
need to show that
$$Y_{s-1}\mathcal{T}_k\cap G^{(s-1)}_s\subsetneq Y_{s-1}\mathcal{T}_k\cap G^{(s-1)}_k\subsetneq Y_{s-1}\mathcal{T}_k.$$
The second inequality is obvious. The first inclusion is strict because of the following reason.
$G^{(s-1)}_s$ is a $\mathcal{G}^{(s-1)}$-factor of $Y_{s-1}\mathcal{T}_k\cap G^{(s-1)}_s$,
therefore $G^{(s-1)}_k$ is not a $\mathcal{G}^{(s-1)}$-factor because it strictly contains
$G^{(s-1)}_s$. On the other hand, $G^{(s-1)}_k$ is a $\mathcal{G}^{(s-1)}$-factor of
$Y_{s-1}\mathcal{T}_k\cap G^{(s-1)}_k$. So the first inclusion is strict.

Case (ii):  $\mathcal{T}_k\cup \{G_s\}$ is not $\mathcal{G}$-nested. In this case, $G^{(s-1)}_s\cap Y_{s-1}\mathcal{T}_k=\emptyset$ by Proposition
\ref{how strata change}. Hence no twisting is needed for the normal bundle.

Case (iii): $\mathcal{T}_k\cup \{G_s\}$ is $\mathcal{G}$-nested but $G_s$ is not strictly contained in $G_k$. If $\mathcal{T}_{k-1}\cup \{G_s\}$ is
not a $\mathcal{G}$-nest, then $G^{(s-1)}_s\cap Y_{s-1}\mathcal{T}_{k-1}=\emptyset$ by Proposition \ref{how strata change}. Hence blowing up along
$G^{(s-1)}_s$ will not affect the normal bundle of $Y_{s-1}\mathcal{T}_{k-1}$, so no twisting is needed. Otherwise, assume $\mathcal{T}_{k-1}\cup
\{G_s\}$ is a $\mathcal{G}$-nest. Both $G_s$ and $G_k$ are minimal in the $\mathcal{G}$-nest $\mathcal{T}_{k-1}\cup\{G_s\}$. Then $G^{(s-1)}_s$ and
$G^{(s-1)}_k$ are minimal in a nest and neither one contains the other, therefore they intersect transversally by the definition of nest. Thus,
$Y_{s-1}\mathcal{T}_{k}\cap G^{(s-1)}_k$ and $Y_{s-1}\mathcal{T}_k\cap G^{(s-1)}_s$, regarded as subvarieties of the ambient space
$Y_{s-1}\mathcal{T}_k$, intersect transversally. Therefore Fact \ref{normal bundle_transversal} applies, and no twisting is needed for the normal
bundle.

\smallskip{Step 3}: Apply the result of Step 2 successively for $s=1,2,\dots,k-1$. The normal bundle
$N_{Y_{k-1}\mathcal{T}_{k-1}}Y_{k-1}\mathcal{T}_{k}$ is isomorphic to
$$\Big{(}g_{k-1,k-1}^*\dots g_{1,k-1}^*
\big{(}N_{Y_0\mathcal{T}_{k-1}}Y_0\mathcal{T}_{k}\big{)}\Big{)}\ts
\big{(}-\sum_{(**)}[G^{(k-1)}_s]|_{Y_{k-1}\mathcal{T}_{k-1}}\big{)}
$$ where the sum is over all $s$ that satisfying condition (**). (Here we have used Lemma \ref{g^-1
divisor}.) Therefore
\begin{align*}&f_k^*g^*_{k,k-1}c(N_{Y_{k-1}\mathcal{T}_{k-1}}Y_{k-1}\mathcal{T}_{k})\\&=
  c\bigg{(}g_{\mathcal{T}}^*\big{(}N_{Y_0\mathcal{T}_{k-1}}Y_0\mathcal{T}_{k}|_{Y_0\mathcal{T}}\big{)}\ts
\mathcal{O}\big{(}-\sum_{(**)}[D_{G_s}]|_{Y_{N}\mathcal{T}_{k-1}}\big{)}\bigg{)}.
\end{align*}

Notice that
$$\big{(}N_{Y_0\mathcal{T}_{k-1}}Y_0\mathcal{T}_{k}\big{)}|_{Y_0\mathcal{T}}=
N_{G_k}(\bigcap_{G_k\subsetneq G\in\mathcal{T}}G)|_{Y_0\mathcal{T}}
$$which is denoted by $N_{G_k}$ by our notation.
(The proof is as follows: Suppose $T_1,...,T_m,T_{m+1},...,T_r$ are the minimal elements of the nest $\mathcal{T}_k$, where the first $m$ elements
contain $G_k$. Then the minimal element of the nest $\mathcal{T}_{k-1}$ are $G_k, T_{m+1},...,T_r$. By the definition of nest, $Y_0\mathcal{T}_k$ is
the transversal intersection $T_1\cap\cdots\cap T_m\cap T_{m+1}\cap\cdots\cap T_r$, and $Y_0\mathcal{T}_{k-1}$ is the transversal intersection
$G_k\cap T_{m+1}\cap\cdots\cap T_r$. Therefore  $$N_{Y_0\mathcal{T}_{k-1}}Y_0\mathcal{T}_k=N_{G_k}(T_1\cap\cdots\cap T_m)|_{Y_0\mathcal{T}_{k-1}}.$$
Since $T_1\cap\cdots\cap T_m=\bigcap_{G_k\subsetneq G\in\mathcal{T}}G$, the conclusion follows immediately.)

Now put everything into Corollary \ref{composition}, we have
\begin{align*}
&\alpha_1\circ\dots\circ\alpha_N\\
&=(j_{\mathcal{T}}\btm g_{\mathcal{T}})_*\prod_{G_k\in \mathcal{T}} \Big{\{} c(g_{\mathcal{T}}^*(N_{G_k})\ts
\mathcal{O}(-\sum_{(*')}[D_{G_s}]|_{Y_{N}\mathcal{T}})
\frac{1}{1+D_{G_k}|_{Y_N\mathcal{T}}}\Big{\}}_{r_k-1-\mu_k}, \\
&\beta_N\circ\dots\circ\beta_1=(g_{\mathcal{T}}\btm j_{\mathcal{T}})_*\prod_{G_k\in
\mathcal{T}}(-D_{G_k})^{\mu_k-1}|_{Y_N\mathcal{T}}.
\end{align*}

Finally, we show that the condition (**) can be replaced by the following condition:

$(\bigstar): G_s\subsetneq G_k$ and $\mathcal{T}\cup\{G_s\}$ is a $\mathcal{G}$-nest.

Indeed, $(\bigstar)$ is stronger than (**). However, for those $G_s$ satisfying (**) but not ($\bigstar$), the divisor
$[D_{G_s}]|_{Y_{N}\mathcal{T}}$ would be trivial because $D_{G_s}\cap Y_N\mathcal{T}=\emptyset$. Therefore, replacing (**) by ($\bigstar$) will not
affect the resulting correspondence.

Hence the proof is complete.
\end{proof}

We write a direct conclusion from Step 3 for later usage:
\begin{cor}\label{normal bundle}
Denote $\pi: G^{(k)}_{k+1}\to G_{k+1}$. Then
$$c(N_{G^{(k)}_{k+1}}Y_k)=c\Big{(}\pi^*N_{(G_{k+1})}Y\ts \sum_{G_{k+1}\supsetneq G\in \mathcal{T}}
(-[D_G])|_{G^{(k)}_{k+1}}\Big{)}.$$
\end{cor}
\begin{proof}
Apply Step 3 to the nest $\mathcal{T}=\{G_{k+1}\}$.
\end{proof}

%%%%%%%%%%%%%%%%%%%%%%%%%%%%%%%%%%%%%%%%%%%%%%%%%%%%%%%%%%%%%%%%%%%%%%%%%%%%%%%%%%

\ss{Fulton-MacPherson configuration spaces} Fix a nonsingular variety $X$ of dimension $d$. The configuration space of $n$ distinct ordered points on
$X$, denoted by $F(X,n)$, can be naturally identified with an open subvariety of the Cartesian product $X^n$:
$$F(X,n):=\{(x_1, x_2,\dots,x_n)\in X^n:\, x_i\neq x_j \hbox{ for }i\neq j\}.$$
In their celebrated paper \cite{FM}, Fulton and MacPherson have discovered an interesting compactification $X[n]$ of the configuration space
$F(X,n)$. The compactification is obtained by replacing the diagonals of $X^n$ by a simple normal crossing divisor. It has many attractive
properties, for example the geometry when $n$ points collide, i.e. the degenerate configuration, can be explicitly described using $X[n]$. $X[n]$ is
closely related to the well known compactification  $\overline{\mathcal{M}}_{0,n}$ of the moduli space of stable rational curves with $n$ marked
points. The reader is referred to the beautiful paper \cite{FM} for the original construction and various applications of the Fulton-MacPherson
configuration space.

The Fulton-MacPherson configuration space $X[n]$ can be realized as a wonderful compactification of an arrangement of subvarieties by taking $Y=X^n$, $\mathcal{G}$ the collection of all diagonals of $X^n$ and therefore the induced arrangement is the set of intersections of diagonals which is called polydiagonals (cf. \cite{Li}).

\sss{Main theorems} First we fix some notations:

(i) Denote $[n]:=\{1,2,\dots,n\}$.  We call two subsets $I, J\subseteq[n]$ \emph{overlapped} if $I\cap J$ is
a nonempty proper subset of $I$ and $J$. For a set $\mathcal{S}$ of subsets of $[n]$, we call
$I$ is compatible with $\mathcal{S}$ (denote by $I\sim \mathcal{S}$) if $I$ does not overlap any
element in $\mathcal{S}$.

A \emph{nest} $\mathcal{S}$ is a set of subsets of $[n]$ such that any two elements $I\neq J\in
\mathcal{S}$ are not overlapped, and all singletons $\{1\},\dots,\{n\}$ are in $\mathcal{S}$.
Notice that the nest defined here, unlike the one defined in \cite{FM}, is allowed to contain
singletons.

Given a nest $\mathcal{S}$, define $\mathcal{S}^\circ=\mathcal{S}\setminus
\big{\{}\{1\},\dots,\{n\}\big{\}}$. In the description of nests by forests below,
$\mathcal{S}^\circ$ correspond to the forest $\mathcal{S}$ cutting of all leaves.

A nest $\mathcal{S}$ naturally corresponds to forest (i.e. a not necessarily connected tree), each node of which is labeled by an element in
$\mathcal{S}$. For example, the following forest corresponds to a nest $\mathcal{S}=\{1,2,3,23,123\}$.
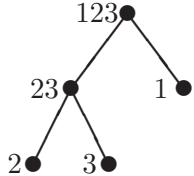
\begin{figure}[h]
\setlength{\unitlength}{.5cm}
\begin{picture}(12,4)\thicklines

\put(0,0){\circle*{0.4}} \put(1,2){\circle*{0.4}} \put(2,0){\circle*{0.4}} \put(2.5,4){\circle*{0.4}} \put(4,2){\circle*{0.4}}
\put(0,0){\line(1,2){1}} \put(1,2){\line(1,-2){1}} \put(2.5,4){\line(-3,-4){1.4}} \put(2.5,4){\line(3,-4){1.4}}

\put(-.5,0){\makebox(0,0){$2$}} \put(0.3,2){\makebox(0,0){$23$}} \put(1.5,0){\makebox(0,0){$3$}} \put(3.4,2){\makebox(0,0){$1$}}
\put(1.7,4){\makebox(0,0){$123$}}
\end{picture}
\end{figure}

Denote by $c(\mathcal{S})$ the number of connected components of the forest, i.e., the number of
maximal elements of $\mathcal{S}$. Denote by $c_I(\mathcal{S})$ (or $c_I$ if no ambiguity arise)
the number of maximal elements of the set $\{J\in\mathcal{S}|J\subsetneq I\}$, i.e. the number of
sons of the node $I$. In the above example, $c(\mathcal{S})=1$, $c_{123}=c_{23}=2$.

(ii) For a subset $I\subseteq [n]$ consisting of at least two elements, define the diagonal
$$\Delta_I:=\{(x_1,\dots,x_n)\in X^n: \, x_i=x_j \hbox{ if } i,j\in I\}.$$
It is shown in \cite{FM} that complement of $F(x,n)$ in the Fulton-MacPherson compactification $X[n]$ is a union of normal crossing nonsingular divisors $D_I$, indexed by subsets $I\subseteq[n]$ with at least two elements. More precisely, $D_I$ is the dominant transform $\widetilde{\Delta}_I$ under the natural morphism $X[n]\to X^n$.

For every nest $\mathcal{S}$, $X(\mathcal{S}):=\cap_{I\in\mathcal{S}}D_I$ is a nonsingular subvariety of $X[n]$. Define $j_\mathcal{S}:
X(\mathcal{S})\hookrightarrow X[n]$ to be the natural inclusion.

Define $\Delta_\mathcal{S}:=\cap_{I\in\mathcal{S}}\Delta_I$. Define
$g_\mathcal{S}:X(\mathcal{S})\to\Delta_\mathcal{S}$ to be the restriction of the morphism $\pi:
X[n]\to X^n$ to the subvariety $X(\mathcal{S})$.

(iii) Let $p_I: X[n]\to X$ be the composition of $\pi:X[n]\to X^n$  with the projection $X^n\to X$ to the $i$-th factor for an
arbitrary $i\in I$. (The choice of $i\in I$ is not essential: indeed, the only place we need $p_I$
is in the formulation of $\alpha_{\mathcal{S}, \mu}$ below, where need the composition
$j_\mathcal{S}^*p_I^*$. By the following diagram
$$\xymatrix{X(\mathcal{S})\ar[r]^{j_\mathcal{S}}\ar[d]^{g_\mathcal{S}} & X[n]\ar[d]^{p_i}\\ \Delta_\mathcal{S}\ar[r]^{q_i} & X}$$
where $i\in I$, we have $j_\mathcal{S}^*p_i^*=g_\mathcal{S}^*q_i^*$, but $q_i$ is independent of
the choice of $i\in I$ since $\Delta_\mathcal{S}\subseteq\Delta_I$, so $j_\mathcal{S}^*p_I^*$ is
independent of the choice of $i\in I$ for $p_I$.)

(iv) For a nest $\mathcal{S}\neq\{\{1\},\dots,\{n\}\}$ (i.e. $\mathcal{S}^\circ\neq\emptyset$), define
$$M_\mathcal{S}:=\big{\{}\underline{\mu}=\{\mu_I\}_{I\in\mathcal{S}^\circ}: 1\le \mu_I\le d(c_I-1)-1,\, \mu_I\in\mathbb{Z} \big{\}}.$$
(recall that $d=\dim X$ and $c_I=c_I(\mathcal{S})$ is defined in (i)) and define
$$||\underline{\mu}||:=\sum_{I\in\mathcal{S}^\circ} \mu_I,\quad \forall \underline{\mu}\in
M_\mathcal{S}.$$
For $\mathcal{S}=\{\{1\},\dots,\{n\}\}$, assume $M_\mathcal{S}=\{\underline{\mu}\}$ with
$\|\underline{\mu}\|=0$.

We will show in the proof of Theorem 4.1 that $M_\mathcal{S}$ is the special case of  $M_\mathcal{T}$ defined in \S \ref{main thm setting} where $Y$ is $X^n$, $\mathcal{G}$ is the set of diagonals of $X^n$ and $\mathcal{T}$ is the set of $\mathcal{G}$-nests.

Define function $\zeta(x):=\sum_{i=0}^d(1+x)^{d-i}c_i(T_X)$.

Define $\alpha_{\mathcal{S}, \underline{\mu}}\in Corr^{-||\underline{\mu}||}(X[n],
\Delta_\mathcal{S})$, $\beta_{\mathcal{S}, \underline{\mu}}\in
Corr^{||\underline{\mu}||}(\Delta_\mathcal{S}, X[n])$, $p_{\mathcal{S}, \underline{\mu}}\in
Corr^0(X[n], X[n])$ as follows,
\begin{align*}
\alpha_{\mathcal{S}, \underline{\mu}}&=(j_\mathcal{S}\btm g_\mathcal{S})_*j_\mathcal{S}^*\bigg{(}\prod_{I\in\mathcal{S}^\circ}
\Big{\{}-p_I^*\zeta(-\sum_{\substack{J\sim\mathcal{S}\\J\supsetneq I}} D_J)^{c_I-1} \frac{1}{1+D_{I}}\Big{\}}_{d(c_I-1)-1-\mu_I}\bigg{)} ,
\\
\beta_{\mathcal{S}, \underline{\mu}}&=(g_\mathcal{S}\btm j_\mathcal{S})_*j_\mathcal{S}^*\bigg{(}\prod_{I\in\mathcal{S}^\circ} D_I^{\mu_I-1}\bigg{)} , \\
p_{\mathcal{S}, \underline{\mu}}&=\beta_{\mathcal{S}, \underline{\mu}}\circ\alpha_{\mathcal{S},
\underline{\mu}}.
\end{align*}
(In the above definition of $\alpha_{\mathcal{S}, \underline{\mu}}$ and $\beta_{\mathcal{S},
\underline{\mu}}$, the products are set to be $1_{X(\mathcal{S})}\in
A^0\big{(}X(\mathcal{S})\big{)}$ if $\mathcal{S}^\circ=\emptyset$.)

\medskip

The following are the main theorems on the Chow groups and Chow motives of Fulton-MacPherson
configuration spaces.

\begin{thm}\label{main FM Chow}
Let $X$ be a nonsingular quasi-projective variety. There is an isomorphism of Chow groups:
$$A^*(X[n])=\bigoplus_{\mathcal{S}}\bigoplus_{\underline{\mu}\in M_{\mathcal{S}}} A^{*-||\underline{\mu}||}(X^{c(\mathcal{S})}) .$$
where $\mathcal{S}$ runs through all nests of $[n]$.
\end{thm}

\begin{thm}\label{main FM motive}
Let $X$ be a nonsingular projective variety. Then there is a canonical isomorphism of Chow motives
$$\bigoplus_{\mathcal{S}}\bigoplus_{\underline{\mu}\in M_{\mathcal{S}}}\alpha_{\mathcal{S}, \underline{\mu}}:
h(X[n])\cong \bigoplus_{\mathcal{S}}\bigoplus_{\underline{\mu}\in
M_{\mathcal{S}}}h(\Delta_\mathcal{S})(\|\underline{\mu}\|)$$ with the inverse
$\sum\limits_{\mathcal{S}}\sum\limits_{\underline{\mu}\in\mathcal{S}}\beta_{\mathcal{S},
\underline{\mu}}$. Equivalently, we have $$h(X[n])\cong
\bigoplus_{\mathcal{S}}\bigoplus_{\underline{\mu}\in M_{\mathcal{S}}}
h(X^{c(\mathcal{S})})(\|\underline{\mu}\|).$$
\end{thm}

\rmk\label{remark symmetry} Observe that the two sets of correspondences $\{\alpha_{\mathcal{S}, \underline{\mu}}\}$, $\{\beta_{\mathcal{S},
\underline{\mu}}\}$ are $\mathfrak{S}_n$-symmetric in the sense that the following diagram commutes for any $\sigma\in
\mathfrak{S}_n$,
$$\xymatrix{h(X[n])\ar[rr]^-{\alpha_{\mathcal{S}, \underline{\mu}}}\ar[d]^\sigma&& h(\Delta_\mathcal{S})(\|\underline{\mu}\|)\ar[rr]^-{\beta_{\mathcal{S}, \underline{\mu}}}\ar[d]^\sigma&& h(X[n])\ar[d]^\sigma \\
h(X[n])\ar[rr]^-{\alpha_{\sigma(\mathcal{S}, \underline{\mu})}}&& h(\Delta_\mathcal{S})(\|\underline{\mu}\|)\ar[rr]^-{\beta_{\sigma(\mathcal{S},
\underline{\mu})}}&& h(X[n])\,.}$$

\begin{proof}[Proof of Theorem \ref{main FM Chow}]

Apply Theorem \ref{main thm wonderful} with the ambient space $Y=X^n$ and the building set
$$\mathcal{G}=\{\Delta_I\}_{I\subseteq [n],|I|\ge 2}$$

First notice that a nest $\mathcal{S}$ of $[n]$ gives a $\mathcal{G}$-nest $\mathcal{T}=\{\Delta_I\}_{I\in\mathcal{S}^\circ}$. Moreover, the inverse
is also true: a $\mathcal{G}$-nest will give a nest of $[n]$. Indeed, given a partition $\Pi=(I_1,\dots,I_t)$ of $[n]$, a $\mathcal{G}$-factor of
$\Delta_{\Pi}$ by definition is a minimal element in $\{G\in\mathcal{G}:\, G\supseteq\Delta_{\Pi}\}$. So $\{\Delta_{I_1},\dots,\Delta_{I_t}\}$ are
all the $\mathcal{G}$-factors of $\Delta_{\Pi}$.  By the definition of $\mathcal{G}$-nest, $\mathcal{T}$ is induced from a flag of strata
$$\Delta_{\Pi_1}\supseteq\Delta_{\Pi_2}\supseteq\dots\supseteq\Delta_{\Pi_t}.$$ Then
$$\Pi_1\ge \Pi_2\ge\dots\ge\Pi_k.$$
(Here $\Pi\ge\Pi'$ means $\Pi$ is a finer partition than $\Pi'$, e.g.$(12,3,4)\ge(123,4)$.) The nest $\mathcal{T}$ is induced by ``taking the union
of all factors of each $\Delta_\Pi$'', which corresponds to ``taking all $I$'s that appear in any of the partition $\Pi_i$''. Since the partitions
are totally ordered, the set of $I$'s forms a nest of $[n]$.

Next we prove the range of $\underline{\mu}$ is as stated. Theorem \ref{main thm wonderful} asserts that
$$1\le\mu_G\le r_G-1.$$ Now $G=\Delta_I$ is a diagonal, so by definition
\begin{align*}r_G&:=\dim(\bigcap_{G\subsetneq
T\in\mathcal{T}}T)-\dim G\\
&=\dim(\bigcap_{I\supsetneq I'\in\mathcal{S}}\Delta_{I'})-\dim\Delta_I\\
&=d(c_I-1).
\end{align*}

Finally, observe that
$$Y_0\mathcal{T}=\bigcap_{G\in\mathcal{T}}G=\bigcap_{I\in\mathcal{S}}\Delta_I=\Delta_\mathcal{S}\cong
X^{c(\mathcal{S})}.$$ Therefore the expected conclusion is implied by Theorem \ref{main thm wonderful}.
\end{proof}

\begin{proof}[Proof of Theorem \ref{main FM motive}]

The statement of the motive decomposition is proved exactly as the above proof.

The correspondences are induced from Theorem \ref{main thm correspondence}.  The improvement of this theorem than Theorem \ref{main thm
correspondence} is: we can say more about the Chern classes appeared in the correspondence $\alpha_{\mathcal{S},\underline{\mu}}$ in Theorem
\ref{main thm correspondence}.

First, given $G=\Delta_I$, let $\Pi=(I_1,\dots,I_{c_I})$ be the partition containing all sons of $I$ in $\mathcal{S}$. We compute the normal bundle
$N_G:=N_{\Delta_I}\Delta_\Pi$. Without loss of generality, assume $I=(12\dots m)$, where $m\le n$.

Denote by $p_i:\Delta_I\to X$ and $q_i:\Delta_\Pi\to X$ the projections induced from the projection of $X^n$ to the $i$-th factor. For each $1\le
i\le c_I$, pick an $a_i\in I_i$.
\begin{align*}
T_{\Delta_I}&=p_1^*T_X\oplus p_{m+1}^*T_X\oplus\dots\oplus p_n^*T_X,\\
T_{\Delta_\Pi}&=q_{a_1}^*T_X\oplus\dots\oplus q_{a_{c_I}}^*T_X\oplus q_{m+1}^*T_X\oplus\dots\oplus q_n^*T_X,\\
T_{\Delta_\Pi}|_{\Delta_I}&=p_1^*T_X\oplus\dots\oplus p_1^*T_X\oplus q_{m+1}^*T_X\oplus\dots\oplus q_n^*T_X,
\end{align*}
therefore, $c(N_G)=p^*_1c(T_X)^{c_I-1}$.

To compute the Chern classes of $N_G$ twisted by a line bundle $L$, we use the Chern root technique. For any vector bundle $N$ on $X$, define the
Chern polynomial as
$$c_y(N):=c_0(N)+c_1(N)y+c_2(N)y^2+\dots.$$
Define $x=c_1(L)$. Recall that the rank of $N_G$ is $r_G=d(c_I-1)$. Then
\begin{align*}
c(N_G\otimes L)&=c_{r_G}(N_G)+c_{r_G-1}(N_G)(1+x)+...+c_0(N_G)(1+x)^{r_G}\\
&=(x+1)^{r_G}c_{\frac{1}{x+1}}(N_G)\\
&=(x+1)^{d(c_I-1)}p_1^*c_{\frac{1}{x+1}}(T_X)^{c_I-1}\\
&=p_1^*[(x+1)^{d}c_{\frac{1}{x+1}}(T_X)]^{c_I-1} =p^*_1\zeta(x)^{c_I-1}.
\end{align*}

Finally, by restricting to $\Delta_\mathcal{S}$ and pulling back to $X(\mathcal{S})$ we get the
expected formula for correspondences $\alpha_{\mathcal{S},\underline{\mu}}$.
\end{proof}

%%%%%%%%%%%%%%%%%%%%%%%%%%%%%%%%%%%%%%%%%%%%%%%%%%%%%%%%%%%%%%%%%%%%%55

\sss{A formula for the generating function of Chow groups and Chow motive of $X[n]$}

In this section, we express the decompositions of the Chow groups (Theorem \ref{main FM Chow})
and the Chow motive (Theorem \ref{main FM motive}) in terms of exponential generating
functions.

Define $[x^i]$ to be the function that picks up the coefficient of $x^i$ from a power series. Define $[\frac{x^it^n}{n!}]$ to be the function that
picks up the coefficient of $\frac{x^it^n}{n!}$ from a power series with two variables $x$ and $t$, i.e.,
$$\Big{[}\frac{x^it^n}{n!}\Big{]}\sum_{j,m}a_{jm}\frac{x^jt^m}{m!}:=a_{in}.$$

The main theorem of this section is the following:
\begin{thm}\label{generating}
Define $f_i(x)$ to be the polynomials whose exponential generating function
$N(x,t)=\sum\limits_{i\ge 1}f_i(x)\frac{t^i}{i!}$ satisfies the identity
$$(1-x)x^dt+(1-x^{d+1})=\exp(x^dN)-x^{d+1}\exp(N).$$
Then for a nonsingular $d$-dimensional quasi-projective variety $X$,
$$A^*(X[n])=\bigoplus_{\substack{1\le k\le n\\i\ge 0}} A^{*-i}(X^k)^{\oplus
[\frac{x^it^n}{n!}]\frac{N^k}{k!}}.$$ Moreover, if $X$ is projective, we have the motive decomposition
\begin{align*}
h(X[n])&=\bigoplus_{\substack{\Pi=(I_1,...,I_k)\\\textrm{partition of } [n]}}
\big{(}h(\Delta_\Pi)(i)\big{)}^{\oplus
[x^i](f_{|I_1|}(x)...f_{|I_k|}(x))}\\
 &=\bigoplus_{\substack{1\le k\le n\\i\ge 0}}\big{(} h(X^k)(i)\big{)}^{\oplus
[\frac{x^it^n}{n!}]\frac{N^k}{k!}}.
\end{align*}
\end{thm}

\rmk
\label{N}
One can write down by hand the first several terms of $N$. Define $\sigma_j=\sum_{i=1}^{dj-1}x^i$
(when $d=1$, define $\sigma_1=0$). Then
\begin{align*}
N&=t+\sigma_1\frac{t^2}{2!}+(\sigma_2+3\sigma_1^2)\frac{t^3}{3!}+(\sigma_3+10\sigma_1\sigma_2+15\sigma_1^3)\frac{t^4}{4!}\\
&+(\sigma_4+15\sigma_1\sigma_3+10\sigma_2^2+105\sigma_1^2\sigma_2+105\sigma_1^4)\frac{t^5}{5!}+....
\end{align*}

\begin{proof}[Proof of Theorem \ref{generating}]

We prove only the statement for motives, since the statement for Chow groups can be proved by exactly the
same method.

By Theorem \ref{main FM motive}, we want to count for any given $i$ and $k$, how many possible
$\mathcal{S}$ and $\underline{\mu}\in\mathcal{S}$ satisfy $c(\mathcal{S})=k$ and
$\|\underline{\mu}\|=i$. First, consider the case when $c(\mathcal{S})=1$, i.e. $\mathcal{S}$ is a
connected forest.

Define
$$f_n(x):=\sum_{\mathcal{S}:c(\mathcal{S})=1}\sum_{\underline{\mu}\in M_\mathcal{S}}x^{||\underline{\mu}||},$$ and
define $f_1(x)=1$.

For a nest $\mathcal{S}$ of $[n]$ with $c(\mathcal{S})=1$, we have
$$\sum_{\underline{\mu}\in M_\mathcal{S}}x^{\|\underline{\mu}\|}=\prod_{I\in\mathcal{S}^\circ}\sigma_{(c_I-1)},$$
i.e., $I$ goes through all non-leaves of $\mathcal{S}$ (if $n=1$, then the sum is set to be
$1$). Since the sons of the root of $\mathcal{S}$ correspond to a partition $\{I_1,\dots,I_k\}$ of
$[n]$, we have following formula for $n\ge 2$,
$$f_n(x)=\sum_{\{I_1,...,I_k\}\\\textrm{partition of } [n]}f_{|I_1|}f_{|I_2|}...f_{|I_k|}\sigma_{k-1} .$$
where $\sigma_{k}=\sum_{i=1}^{dk-1}x^i$ for $k>0$, and $\sigma_0=0$. Since the equality does not
hold for $n=1$ where $f_1(x)=1$ but the right side is $0$, so one define
$$\tilde{f}_n(x)=\left\{
             \begin{array}{ll}
               f_n(x), & \hbox{if $n>1$;} \\
               0, & \hbox{if $n=1$.}
             \end{array}
           \right.
$$ Then the following holds for any $n\ge 1$:
$$\tilde{f}_n(x)=\sum_{\{I_1,...,I_k\}\\\textrm{partition of } [n]}f_{|I_1|}f_{|I_2|}...f_{|I_k|}\sigma_{k-1} .$$

Recall the Compositional Formula of exponential generating functions (cf. \cite{Stanley}, Theorem
5.1.4), which asserts that if an equation as above holds, then
$$E_{\tilde{f}}(t)=E_\sigma(E_f(t)),$$
where
\begin{align*}
E_{\tilde{f}}(t)&=1+\tilde{f}_1t+\tilde{f}_2t^2/2!+\tilde{f}_3t^3/3!+\dots\\
E_{\sigma}(t)&=1+\sigma_0t+\sigma_1t^2/2!+\sigma_2t^3/3!+\dots\\
E_{f}(t)&=\qquad f_1t+f_2t^2/2!+f_3t^3/3!+\dots
\end{align*}

By the definition of $\tilde{f}$, $E_{\tilde{f}}=E_f-t+1$. Denote $N=E_{f}$, one has
$$N-t+1=E_g(N),$$
Standard Computation shows
$$E_g(N)=1+N+\frac{1}{x-1}\big{[}\frac{1}{x^d}(e^{x^dN}-1)-xe^N+x\big{]}.$$
Therefore $$(1-x)x^dt+(1-x^{d+1})=\exp(x^dN)-x^{d+1}\exp(N).$$

Now consider the case when $c(\mathcal{S})$ is not necessarily $1$, i.e., the forest $\mathcal{S}$
is not necessarily connected. For a partition $\Pi=\{I_1,...,I_k\}$ of $[n]$, the number of times
that $h(\Delta_\Pi)(i)$ appears in the decomposition of $h(X[n])$ is equal to
$[x^k](f_{|I_1|}(x)...f_{|I_k|}(x))$, the coefficient of $x^k$ in the product. Denote by $a_{k,i}$
the sum of these numbers for all partitions with $k$ blocks. Then $a_{k,i}$ is the number of times
that $h(X^k)(i)$ appears in the decomposition of $H(X[n])$.

Define $$F_n(y)=\sum_{\{I_1,...,I_k\}\\\textrm{partition of } [n]}f_{|I_1|}f_{|I_2|}...f_{|I_k|}y^k
.$$ Then the coefficient $[y^k]F_n(y)=\sum a_{k,i}x^i$. Use the Compositional Formula again,
$$F_n=[\frac{t^n}{n!}]\exp(yN).$$
Therefore
\begin{align*}
[y^k]F_n(y)&=[y^k][\frac{t^n}{n!}]\exp(yN)\\
&=[\frac{t^n}{n!}][y^k]\exp(yN)\\
&=[\frac{t^n}{n!}]\frac{N^k}{k!}.
\end{align*}
This yields the formula for the decomposition of the Chow motive $h(X[n])$.
\end{proof}

%%%%%%%%%%%%%%%%%%%%%%%%%%%%%%%%%%%%%%%%%%%%%%%%%%%%%%%%%%%%%%%%%5

\sss{Description of $X[n]$ for small $n$}\label{small n}

In this section we explain the previous Theorems (\ref{main FM Chow}, \ref{main FM motive}, and
\ref{generating}) about Fulton-MacPherson configuration space $X[n]$ for small $n=2,3,4$.

For unification of expression, assume $d>1$ in the following examples. (The case $d=1$ is simpler but the expression needs to be modified.)

\smallskip
\noindent\textbf{Example $n=2$}. The morphism $\pi: X[2]\to X^2$ is a blow-up along the diagonal $\Delta_{12}$. Theorem \ref{generating} asserts
\begin{equation}\label{eq:X[2]}
h(X[2])\cong h(X^2)\oplus\bigoplus_{i=1}^{d-1} h(\Delta_{12})(i) \cong h(X^2)\oplus\bigoplus_{i=1}^{d-1} h(X)(i).
\end{equation}
There are 2 possible nests: $\mathcal{S}=\{1,2\}$ and $\mathcal{S}=\{1,2,12\}$. Theorem \ref{main FM motive} asserts the following.

For the first nest, $M_\mathcal{S}$ contains only one element $\underline{\mu}$ with $\|\underline{\mu}\|=0$. Therefore $\alpha=\Gamma_\pi$,
$\beta=\Gamma^t_\pi$, $p=\Gamma^t_\pi\circ\Gamma_\pi$. They give the first direct summand in the decomposition (\ref{eq:X[2]}).

For the second nest, $\mathcal{S}^\circ=\{12\}$, $1\le\mu_{12}\le d-1$, so there are $d-1$ direct summands
for this nest. Denote $j: D_{12}\hookrightarrow X[2]$, $g: D_{12}\to\Delta_{12}$ as the natural map, we have
\begin{align*}
\alpha_{\mathcal{S}, \underline{\mu}}&=-(j\btm g)_*j^*\big{(} \sum_{i=0}^{d-1-\mu_{12}}p_1^*c_i(T_X)
(-D_{12})^{d-1-\mu_{12}-i} \big{)},\\
\beta_{\mathcal{S}, \underline{\mu}}&=(g\btm j)_*j^*\big{(}D^{\mu_{12}-1}{)} , \\
p_{\mathcal{S}, \underline{\mu}}&=\beta_{\mathcal{S}, \underline{\mu}}\circ\alpha_{\mathcal{S},
\underline{\mu}} .
\end{align*}
They give the direct summand $h(\Delta_{12})(\mu_{12})$ in the decomposition (\ref{eq:X[2]}).
\begin{figure}[b]
\label{fig1}
% Example of a chain of partitions and its leveled tree
\hskip10mm\xy
    (88,-8)*{\mycaption{7in}{
            $X[3]$ by the symmetric construction.}};
% Annotate partitions
 %   (-7,14.7)*{\scriptstyle \pi_1};
 %   (-7, 9.7)*{\scriptstyle \pi_2};
 %   (-7, 4.7)*{\scriptstyle \pi_3};
%   (-7,-0.3)*{\scriptstyle \top};
% Ugly tree
    (24,25)="root"      ;
    (18,20)="a1"    ;
    (29,20)="a2"        *{\scriptstyle 23};
    (12,15)="b1"    ;
    (23,15)="b2"        *{\scriptstyle 13};
    (35,15)="b3"       *{\scriptstyle 23};
    (6,10)="c1"     ;
    (17,10)="c2"      *{\scriptstyle 12};
    (29,10)="c3"    *{\scriptstyle 13};
    (41,10)="c4"        *{\scriptstyle 23};
    (2,5)="d1"        ;
    (8,5)="d2"        *{\scriptstyle 123};
    (14,5)="d3"       *{\scriptstyle 12};
    (20,5)="d4"     *{\scriptstyle 123};
    (26,5)="d5"      *{\scriptstyle 13};
    (32,5)="d6"       *{\scriptstyle 123};
    (38,5)="d7"      *{\scriptstyle 23};
    (44,5)="d8"      *{\scriptstyle 123};
    (2,0)="e1"        *{\textrm{\textcircled{\raisebox{-1pt}{1}}}};
    (8,0)="e2"        *{\textrm{\textcircled{\raisebox{-1pt}{5}}}};
    (14,0)="e3"       *{\textrm{\textcircled{\raisebox{-1pt}{2}}}};
    (20,0)="e4"       *{\textrm{\textcircled{\raisebox{-1pt}{6}}}};
    (26,0)="e5"       *{\textrm{\textcircled{\raisebox{-1pt}{3}}}};
    (32,0)="e6"       *{\textrm{\textcircled{\raisebox{-1pt}{7}}}};
    (38,0)="e7"       *{\textrm{\textcircled{\raisebox{-1pt}{4}}}};
    (44,0)="e8"       *{\textrm{\textcircled{\raisebox{-1pt}{8}}}};
    "root";         "a1"**\dir{-};
    "root";         "a2"+(-0.3,1.2)**\dir{-};
    "a1";           "b1"**\dir{-};
    "a1";           "b2"+(0,1.2)**\dir{-};
    "a2"-(-2,1.2);  "b3"+(-0.5,1.2)**\dir{-};
    "b1";           "c1"**\dir{-};
    "b1";           "c2"+(-0.5,1.2)**\dir{-};
    "b2"-(-2,1.2);  "c3"+(-0.5,1.2)**\dir{-};
    "b3"-(-2,1.2);  "c4"+(-.8,1.2)**\dir{-};
    "c1";           "d1"+(0,1)**\dir{-};
    "c1";           "d2"+(-0.5,1.2)**\dir{-};
    "c2"-(0.5,1.2);           "d3"+(0,1.2)**\dir{-};
    "c2"-(0,1.2);           "d4"+(-0.5,1.2)**\dir{-};
    "c3"-(0.5,1.2);  "d5"+(0,1.2)**\dir{-};
    "c3"-(0,1.2);  "d6"+(-0.5,1.2)**\dir{-};
    "c4"-(0.5,1.2);  "d7"+(0,1.2)**\dir{-};
    "c4"-(0,1.2);  "d8"+(-0.5,1.2)**\dir{-};
% Indicate levels
    (54,25)*\txt{\tiny level 4};
    (54,20)*\txt{\tiny level 3};
    (54,15)*\txt{\tiny level 2};
    (54,10)*\txt{\tiny level 1};
    (54, 5)*\txt{\tiny level 0};
%
%
%
% Different nests
%1
    (70,0)="m";
    "m"+(3,29)  *{\textrm{\textcircled{\raisebox{-1pt}{1}}}};
    "m"+(0,25)  *{\bullet};
    "m"+(3,25)  *{\bullet};
    "m"+(6,25)  *{\bullet};
    "m"+(0,22)  *{\scriptstyle 1};
    "m"+(3,22)  *{\scriptstyle 2};
    "m"+(6,22)  *{\scriptstyle 3};
%2
    "m"+(15,0)="m";
    "m"+(4,29)  *{\textrm{\textcircled{\raisebox{-1pt}{2}}}};
    "m"+(2,25)="12"  *{\bullet};
    "m"+(8,25)  *{\bullet};
    "m"+(0,22)="1"  *{\bullet};
    "m"+(4,22)="2"  *{\bullet};
    "m"+(0,20)  *{\scriptstyle 1};
    "m"+(4,20)  *{\scriptstyle 2};
    "m"+(8,23)  *{\scriptstyle 3};
    "12";         "1"**\dir{-};
    "12";         "2"**\dir{-};
%3
    "m"+(15,0)="m";
    "m"+(4,29)   *{\textrm{\textcircled{\raisebox{-1pt}{3}}}};
    "m"+(2,25)="12"  *{\bullet};
    "m"+(8,25)  *{\bullet};
    "m"+(0,22)="1"  *{\bullet};
    "m"+(4,22)="2"  *{\bullet};
    "m"+(0,20)  *{\scriptstyle 1};
    "m"+(4,20)  *{\scriptstyle 3};
    "m"+(8,23)  *{\scriptstyle 2};
    "12";         "1"**\dir{-};
    "12";         "2"**\dir{-};
%4
    "m"+(15,0)="m";
    "m"+(4,29)   *{\textrm{\textcircled{\raisebox{-1pt}{4}}}};
    "m"+(2,25)="12"  *{\bullet};
    "m"+(8,25)  *{\bullet};
    "m"+(0,22)="1"  *{\bullet};
    "m"+(4,22)="2"  *{\bullet};
    "m"+(0,20)  *{\scriptstyle 2};
    "m"+(4,20)  *{\scriptstyle 3};
    "m"+(8,23)  *{\scriptstyle 1};
    "12";         "1"**\dir{-};
    "12";         "2"**\dir{-};
%5
    "m"+(-45,-20)="m";
    "m"+(4,29)   *{\textrm{\textcircled{\raisebox{-1pt}{5}}}};
    "m"+(4,25)="123"  *{\bullet};
    "m"+(0,22)="1"  *{\bullet};
    "m"+(4,22)="2"  *{\bullet};
    "m"+(8,22)="3"  *{\bullet};
    "m"+(0,20)  *{\scriptstyle 1};
    "m"+(4,20)  *{\scriptstyle 2};
    "m"+(8,20)  *{\scriptstyle 3};
    "123";         "1"**\dir{-};
    "123";         "2"**\dir{-};
    "123";         "3"**\dir{-};
%6
    "m"+(15,0)="m";
    "m"+(5,29)   *{\textrm{\textcircled{\raisebox{-1pt}{6}}}};
    "m"+(6,25)="123"  *{\bullet};
    "m"+(3,22)="12"  *{\bullet};
    "m"+(9,22)="3"  *{\bullet};
    "m"+(1,19)="1"  *{\bullet};
    "m"+(5,19)="2"  *{\bullet};
    "m"+(1,17)  *{\scriptstyle 1};
    "m"+(5,17)  *{\scriptstyle 2};
    "m"+(9,20)  *{\scriptstyle 3};
    "123";         "12"**\dir{-};
    "123";         "3"**\dir{-};
    "12";         "1"**\dir{-};
    "12";         "2"**\dir{-};
%7
    "m"+(15,0)="m";
    "m"+(5,29)  *{\textrm{\textcircled{\raisebox{-1pt}{7}}}};
    "m"+(6,25)="123"  *{\bullet};
    "m"+(3,22)="12"  *{\bullet};
    "m"+(9,22)="3"  *{\bullet};
    "m"+(1,19)="1"  *{\bullet};
    "m"+(5,19)="2"  *{\bullet};
    "m"+(1,17)  *{\scriptstyle 1};
    "m"+(5,17)  *{\scriptstyle 3};
    "m"+(9,20)  *{\scriptstyle 2};
    "123";         "12"**\dir{-};
    "123";         "3"**\dir{-};
    "12";         "1"**\dir{-};
    "12";         "2"**\dir{-};
%8
    "m"+(15,0)="m";
    "m"+(5,29)   *{\textrm{\textcircled{\raisebox{-1pt}{8}}}};
    "m"+(6,25)="123"  *{\bullet};
    "m"+(3,22)="12"  *{\bullet};
    "m"+(9,22)="3"  *{\bullet};
    "m"+(1,19)="1"  *{\bullet};
    "m"+(5,19)="2"  *{\bullet};
    "m"+(1,17)  *{\scriptstyle 2};
    "m"+(5,17)  *{\scriptstyle 3};
    "m"+(9,20)  *{\scriptstyle 1};
    "123";         "12"**\dir{-};
    "123";         "3"**\dir{-};
    "12";         "1"**\dir{-};
    "12";         "2"**\dir{-};
\endxy
\end{figure}

\smallskip
\noindent\textbf{Example $n=3$}. Apply Theorem \ref{generating},
\begin{align*}
h(X[3])&\cong h(X^3)\oplus\bigoplus_{i=1}^{d-1} h(\Delta_{12})(i)\oplus
\bigoplus_{i=1}^{d-1}h(\Delta_{13})(i)\oplus\bigoplus_{i=1}^{d-1}h(\Delta_{23})(i)\\
&\quad\oplus\bigoplus_{i=1}^{2d-1}\big{(}h(\Delta_{123})(i)\big{)}^{\oplus min\{3i-2, 6d-3i-2\}}\\
&\cong h(X^3)\oplus\bigoplus_{i=1}^{d-1}\big{(}h(X^2)(i)\big{)}^{\oplus
3}\oplus\bigoplus_{i=1}^{2d-1}\big{(}h(X)(i)\big{)}^{\oplus min\{3i-2, 6d-3i-2\}}
\end{align*}

Now we write out all the correspondences that give the decomposition of motives. There are $8$
possible nests, correspond to $8$ trees (see the right side of Figure \ref{fig1}).

The tree on the left side of Figure \ref{fig1} helps us to understand the relation between
subvarieties of different $Y_i$'s (i.e. at different levels): each node with label $I$ at level $k$
correspond to the subvariety $Y_kI:=(\Delta_I)^{(k)}$ in $Y_k$. The node at level $k$ without label
correspond to $Y_k$. For example, the root at level 4 corresponds to $Y_4$, its two successors
correspond to $Y_3$ and $Y_3(23)$, and the relation is that $Y_4$ is the blow-up of $Y_3$ along
$Y_3(23)$.

We list below those correspondences $\alpha, \beta, p $ for the $8$ trees:

\textcircled{\footnotesize{1}} gives $\alpha=\Gamma_\pi, \beta=\Gamma^t_\pi, p=\Gamma^t_\pi\circ\Gamma_\pi$.

\textcircled{\footnotesize{2}} (and \textcircled{\footnotesize{3}}, \textcircled{\footnotesize{4}} are
similar) gives
\begin{align*}
\alpha_{\mathcal{S}, \underline{\mu}}&=(j_\mathcal{S}\btm g_\mathcal{S})_*j_\mathcal{S}^*\big{(}\{-p_1^*\zeta(-D_{123})
\frac{1}{1+D_{12}}\}_{d-1-\mu_{12}}\big{)} ,\\
\beta_{\mathcal{S}, \underline{\mu}}&=(g_\mathcal{S}\btm j_\mathcal{S})_*j_\mathcal{S}^*\big{(}D_{12}^{\mu_{12}-1}\big{)}.
\end{align*}
where $X(\mathcal{S})=D_{12}$, $1\le\mu_{12}\le d-1$.

\textcircled{\footnotesize{5}} gives
\begin{align*}
\alpha_{\mathcal{S}, \underline{\mu}}&=(j_\mathcal{S}\btm g_\mathcal{S})_*j_\mathcal{S}^*\big{(}\{-p_1^*\zeta(\O)^2
\frac{1}{1+D_{123}}\}_{2d-1-\mu_{123}}\big{)} ,\\
\beta_{\mathcal{S}, \underline{\mu}}&=(g_\mathcal{S}\btm j_\mathcal{S})_*j_\mathcal{S}^*\big{(}D_{123}^{\mu_{123}-1}\big{)}.
\end{align*}
where $X(\mathcal{S})=D_{123}$, $1\le\mu_{123}\le 2d-1$.

\textcircled{\footnotesize{6}} (and \textcircled{\footnotesize{7}}, \textcircled{\footnotesize{8}} are
similar) gives
\begin{align*}
\alpha_{\mathcal{S}, \underline{\mu}}=&(j_\mathcal{S}\btm
g_\mathcal{S})_*j_\mathcal{S}^*\\
&\big{(}\{p_1^*\zeta(-D_{123}) \frac{1}{1+D_{12}}\}_{d-1-\mu_{12}}\{p_1^*\zeta(\O)
\frac{1}{1+D_{123}}\}_{d-1-\mu_{123}}\big{)} ,\\
\beta_{\mathcal{S}, \underline{\mu}}=&(g_\mathcal{S}\btm j_\mathcal{S})_*j_\mathcal{S}^*\big{(}D_{12}^{\mu_{12}-1}D_{123}^{\mu_{123}-1}\big{)}.
\end{align*}
where $X(\mathcal{S})=D_{12}\cap D_{123}$, $1\le\mu_{12},\mu_{123}\le d-1$.

\rmk If we use Fulton and MacPherson's nonsymmetric construction of $X[3]$, we would get another set of
correspondences which also gives a decomposition of the motive $h(X[n])$. This set of
correspondences turns out to be different than the ones given above: a straightforward calculation
shows that, by the nonsymmetric construction of $X[3]$, the correspondence that gives the direct
summand $h(\Delta_{12})(\mu_{12})$ is
\begin{align*}
&\alpha:\quad h(X[3])\to h(\Delta_{12})(\mu_{12}),\\
&\alpha=(j_{12}\btm g_{12})_*j_{12}^*\big{(}\{p_1^*\zeta(\O) \frac{1}{1+D_{12}}\}_{d-1-\mu_{12}}\big{)}.
\end{align*}
where $j_{12}:D_{12}\hookrightarrow X[3]$ and $g_{12}:D_{12}\to \Delta_{12}$ are the natural
morphisms. However, the correspondence giving the direct summand $h(\Delta_{13})(\mu_{13})$ is
\begin{align*}
&\alpha':\quad h(X[3])\to h(\Delta_{13})\ts\mathbb{L}^{\mu_{13}},\\
&\alpha'=(j_{13}\btm g_{13})_*j_{13}^*\big{(}\{p_1^*\zeta(-D_{123}) \frac{1}{1+D_{13}}\}_{d-1-\mu_{13}}\big{)}.
\end{align*}
where $j_{13}:D_{13}\hookrightarrow X[3], g_{13}:D_{13}\to \Delta_{13}$ are the natural morphisms. Notice that $\alpha$ and $\alpha'$ are not of
similar forms (Compare $\zeta(\O)$ with $\zeta(-D_{123})$). Therefore the non-symmetry of the construction of $X[3]$ induces the non-symmetry of
correspondences. Actually, this is one reason why we choose the symmetric construction of $X[n]$ (cf. Remark \ref{remark symmetry}).

\smallskip
\noindent\textbf{Example $n=4$}.  we only look at one nest $\mathcal{S}$:
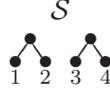
\begin{figure}[h]
\hskip10mm\xy
    (0,0)="m";
    "m"+(6,14) *{\mathcal{S}};
    "m"+(2,10)="12"      *{\bullet};
    "m"+(10,10)="34"      *{\bullet};
    "m"+(0,7)="1"        *{\bullet};
    "m"+(4,7)="2"        *{\bullet};
    "m"+(8,7)="3"        *{\bullet};
    "m"+(12,7)="4"       *{\bullet};
    "m"+(0,5)       *{\scriptstyle 1};
    "m"+(4,5)       *{\scriptstyle 2};
    "m"+(8,5)       *{\scriptstyle 3};
    "m"+(12,5)       *{\scriptstyle 4};
  "12";         "1"**\dir{-};
  "12";         "2"**\dir{-};
  "34";         "3"**\dir{-};
  "34";         "4"**\dir{-};
\endxy
\end{figure}

We have $X(\mathcal{S})=D_{12}\cap D_{34}$, $1\le \mu_{12},\mu_{34}\le d-1$ and
\begin{align*}
\alpha_{\mathcal{S}, \underline{\mu}}=&(j_\mathcal{S}\btm
g_\mathcal{S})_*j_\mathcal{S}^*\\
&\big{(}\{p_1^*\zeta(-D_{1234}) \frac{1}{1+D_{12}}\}_{d-1-\mu_{12}}\{p_3^*\zeta(-D_{1234})
\frac{1}{1+D_{34}}\}_{d-1-\mu_{34}}\big{)} ,\\
\beta_{\mathcal{S}, \underline{\mu}}=&(g_\mathcal{S}\btm j_\mathcal{S})_*j_\mathcal{S}^*\big{(}D_{12}^{\mu_{12}-1}D_{34}^{\mu_{34}-1}\big{)}.
\end{align*}
Since $\Delta_{12}$ and $\Delta_{34}$ would not be disjoint in the procedure of blow-ups, a priori we have to make a choice of order that whether
blow up along (the strict transform of) $\Delta_{12}$ first, or along (the strict transform of) $\Delta_{34}$ first. Although we have to choose
(non-canonically) an order to compute the correspondences, it turns out that the correspondences (hence projectors) which give the motive
decomposition in Theorem \ref{main FM motive} are actually independent of the choice, therefore ``canonical". This independence is a special case of
Remark \ref{remark symmetry}: for $\sigma=(13)(24)\in \mathbb{S}_4$, the above correspondences is invariant under the action induced by $\sigma$.

An application of Theorem \ref{generating} is: we can compute the rank of $A(X[n])$ (as an abelian group) once given the ranks of $A(X^k)$ for all
$1\le k\le n$  (assuming that the ranks of $A(X^k)$'s are finite).

Let us take $\mathbb{P}^d[5]$ for example. Since the rank of $A((\mathbb{P}^d)^k)$ is $(d+1)^k$,
Theorem \ref{generating} implies that the rank of $A(\mathbb{P}^d[5])$ is
$$\sum_{1\le k\le 5}(d+1)^k\Big{(}[\frac{t^5}{t!}]\big{(}\frac{N^k}{k!}|_{x=1}\big{)}\Big{)}.$$
By Remark \ref{N}, we can compute the following
\begin{align*}
\frac{N^2}{2!}&=\frac{t^2}{2!}+3\sigma_1\frac{t^3}{3!}+(15\sigma_1^2+4\sigma_2)\frac{t^4}{4!}
+(105\sigma_1^3+60\sigma_1\sigma_2+5\sigma_3)\frac{t^5}{5!}+.... \\
\frac{N^3}{3!}&=\frac{t^3}{3!}+6\sigma_1\frac{t^4}{4!} +(45\sigma_1^2+10\sigma_2)\frac{t^5}{5!}+....\\
\frac{N^4}{4!}&=\frac{t^4}{4!}+10\sigma_1\frac{t^5}{5!}+....\\
\frac{N^5}{5!}&=\frac{t^5}{5!}+....
\end{align*}
Now plug in $x=1$, we have $\sigma_j=dj-1$. The above sum is a polynomial of $d$ as follows
\begin{align*}
&(d+1)^5 +(d+1)^4 10\sigma_1+(d+1)^3(45\sigma_1^2+10\sigma_2)\\
&+(d+1)^2(105\sigma_1^3+60\sigma_1\sigma_2+5\sigma_3)\\
&+(d+1)(\sigma_4+15\sigma_1\sigma_3+10\sigma_2^2+105\sigma_1^2\sigma_2+105\sigma_1^4).
\end{align*} In particular, the rank of $A(\mathbb{P}^1[5])$ is $178$, the rank of $A(\mathbb{P}^2[5])$ is $7644$.

\rmk For the example $X=\mathbb{P}^d$, since $X[n]$ has an affine cell decomposition, the rank of the Chow group $A_k(X[n])$ coincides with the
$2k$-th Betti number of $X[n]$. Therefore we could also get the above rank by the Poincar\'e polynomial of $X[n]$ computed in \cite{FM}. However, the
rank of $A(X[n])$ for a general variety $X$ is not implied by the Poincar\'e polynomial of $X[n]$.

\ss{Chow motives of $X[n]/\mathfrak{S}_n$}

It is proved in \cite{FM} that the isotropy group of any point in $X[n]$ is a solvable group. It is natural to consider the quotient space
$X[n]/\mathfrak{S}_n$. In this section, we compute its Chow motive in terms of the Chow motives of the Cartesian products of symmetric products of
$X$.

The base field is of characteristic 0 throughout this section.

\begin{lem}\label{quot lem 1}
  Suppose a finite group $G$ acts on a nonsingular projective variety $Y$. If $p_1,\dots, p_k$ are
  orthogonal projectors of $Y$ that

  i) $\sigma p_i=p_i\sigma$, $\forall 1\le i\le k, \forall\sigma\in G$.

  ii) $p_1+p_2+\cdots p_k=\Delta_Y$.

  Then $\ave\Delta_Y=\sum \ave \circ p_i$ where $\ave \circ  p_1, \dots, \ave \circ  p_k$ are orthogonal projectors.
  Consequently, $h(Y)=\oplus(Y, \ave \circ  p_i)$.
\end{lem}

\begin{proof} Since
  \begin{align*}
    (\ave\, p_i) (\ave\, p_j)&=\Big{(}\dfrac{1}{|G|}\sum_\sigma \sigma p_i\Big{)}\Big{(}\dfrac{1}{|G|}\sum_\tau
    \tau
    p_j\Big{)}\\
    &=\dfrac{1}{|G|^2}\sum_{\sigma, \tau} \sigma\tau p_i p_j
    =\dfrac{1}{|G|}\sum_{\sigma} \sigma \delta_{ij}p_i
    =\delta_{ij}(\ave\, p_j).
  \end{align*} Then the lemma follows.
\end{proof}

\begin{lem}\label{quot lem 2}
  Suppose $Y, Z$ are nonsingular (not necessary connected) projective varieties with finite group $G$
  actions. Suppose that $\alpha\in Corr^{-m}(Y, Z)$ has an inverse $\beta\in Corr^m(Z,Y)$, and $\alpha$ gives an isomorphism of Chow
  motives $$(Y, p)\cong h(Z)(m)$$ where $p=\beta\alpha$, and $\alpha\sigma=\sigma\alpha$, $\beta\sigma=\sigma\beta$,  $\forall
  \sigma\in G$. Then $$(Y, \ave\circ p)\cong h(Z/G)(m).$$
\end{lem}

\begin{proof}
  Similar to the proof of Lemma \ref{quot lem 1}, we have $(\ave\, p)^2=\ave\, p$ and the following
  commutative diagram
  $$\xymatrix{Y\ar[r]^-{\ave\, \alpha}\ar[d]_{\ave\, p} & Z\ar[r]^-{\ave\,\beta}\ar[d]^{\ave\, \Delta_Z}& Y\ar[d]^{\ave\, p} \\
  Y\ar[r]_-{\ave\, \alpha} & Z\ar[r]_-{\ave\,\beta} & Y.}$$
  Therefore, $(Y, \ave\, p)\cong (Z, \ave\,\Delta_Z)(m)\cong h(Z/G)(m).$
\end{proof}

Now we consider the quotient variety $X[n]/\mathfrak{S}_n$. For convenience, define $G:=\mathfrak{S}_n$. There is a natural action of $G$ on the set
$\{(\mathcal{S},\underline{\mu})\}$ where $\mathcal{S}$ are nests and $\underline{\mu}\in M_\mathcal{S}$.  Define the subgroup
$G_{\mathcal{S},\underline{\mu}}$ of $G$ as
$$G_{\mathcal{S},\underline{\mu}}=\{\sigma\in\frak{S}_n:\; \sigma(\mathcal{S},\underline{\mu})=(\mathcal{S},\underline{\mu})\}.$$
Define $\overline{(\mathcal{S},\underline{\mu})}$ to be the class of $G$-orbit $G\cdot
(\mathcal{S},\underline{\mu})$. Then
$$\Delta_Y=\sum_{\mathcal{S},\underline{\mu}}p_{\mathcal{S},\underline{\mu}}
=\sum_{\overline{(\mathcal{S},\underline{\mu})}}\quad\sum_{\overline{\sigma}\in
G/G_{\mathcal{S},\underline{\mu}}}p_{\sigma(\mathcal{S},\underline{\mu})}.$$ Since $\{\alpha_{\mathcal{S}, \underline{\mu}}\}$,
$\{\beta_{\mathcal{S}, \underline{\mu}}\}$ are $\mathfrak{S}_n$-symmetric (cf. the Remark after Theorem \ref{main FM motive}), it is easy to check
that $\sum_{\overline{\sigma}\in G/G_{\mathcal{S},\underline{\mu}}}p_{\sigma(\mathcal{S},\underline{\mu})}$ commutes with every $\tau\in G$. By Lemma
\ref{quot lem 1},
$$h(X[n]/G)\cong(Y, \ave\circ\Delta_Y)\cong\bigoplus_{\overline{(\mathcal{S},\underline{\mu})}}\Big{(}Y, \ave\circ\sum_{\overline{\sigma}\in
G/G_{\mathcal{S},\underline{\mu}}}p_{\sigma(\mathcal{S},\underline{\mu})}\Big{)}.$$ Since
$$\Big{(}Y, \sum_{\overline{\sigma}\in
G/G_{\mathcal{S},\underline{\mu}}}p_{\sigma(\mathcal{S},\underline{\mu})}\Big{)}\cong\Big{(}
\bigsqcup_{\overline{\sigma}\in
G/G_{\mathcal{S},\underline{\mu}}}\Delta_{\sigma(\mathcal{S})}\Big{)}(||\underline{\mu}||),$$ by
Lemma \ref{quot lem 2} we have
\begin{align*}
\Big{(}Y, \ave\circ\sum_{\overline{\sigma}\in
G/G_{\mathcal{S},\underline{\mu}}}p_{\sigma(\mathcal{S},\underline{\mu})}\Big{)}&\cong
h\bigg{(}\Big{(}\bigsqcup_{\overline{\sigma}\in
G/G_{\mathcal{S},\underline{\mu}}}\Delta_{\sigma(\mathcal{S})}\Big{)}{/G}\bigg{)}(||\underline{\mu}||)\\
&\cong h(\Delta_{\mathcal{S}}/G_{\mathcal{S},\mu})(||\underline{\mu}||).
\end{align*}

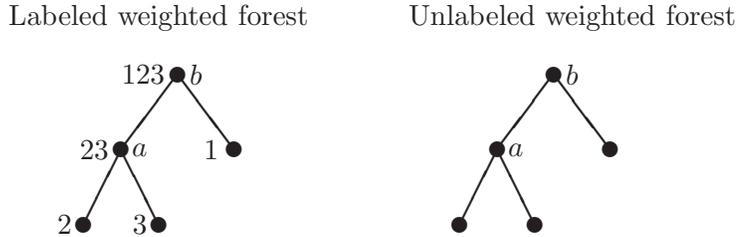
\begin{figure}[h]
\setlength{\unitlength}{.5cm}
\begin{picture}(12,6)\thicklines

\put(0,0){\circle*{0.4}} \put(1,2){\circle*{0.4}} \put(2,0){\circle*{0.4}} \put(2.5,4){\circle*{0.4}} \put(4,2){\circle*{0.4}}
\put(0,0){\line(1,2){1}} \put(1,2){\line(1,-2){1}} \put(2.5,4){\line(-3,-4){1.4}} \put(2.5,4){\line(3,-4){1.4}}

\put(-.5,0){\makebox(0,0){$2$}} \put(0.3,2){\makebox(0,0){$23$}} \put(1.5,0){\makebox(0,0){$3$}} \put(3.4,2){\makebox(0,0){$1$}}
\put(1.6,4){\makebox(0,0){$123$}}

\put(1.5,2){\makebox(0,0){$a$}} 
\put(3,4){\makebox(0,0){$b$}}
\put(2,5.5){\makebox(0,0){Labeled weighted forest}}

\put(10,0){\circle*{0.4}} \put(11,2){\circle*{0.4}} \put(12,0){\circle*{0.4}} \put(12.5,4){\circle*{0.4}} \put(14,2){\circle*{0.4}}
\put(10,0){\line(1,2){1}} \put(11,2){\line(1,-2){1}} \put(12.5,4){\line(-3,-4){1.4}} \put(12.5,4){\line(3,-4){1.4}}

\put(11.5,2){\makebox(0,0){$a$}} 
\put(13,4){\makebox(0,0){$b$}}
\put(13,5.5){\makebox(0,0){Unlabeled weighted forest}}
\end{picture}
 \caption{\label{figure:weighted forest} Labeled and unlabeled weighted forests.}
\end{figure}

The space $\Delta_{\mathcal{S}}/G_{\mathcal{S},\underline{\mu}}$ can be described as follows. Each
$(\mathcal{S},\underline{\mu})$ corresponds to a labeled ``weighted'' forest, the correspondence is given 
by attaching  an integer $\mu_I$ to each non-leaf node $I$ of the labeled forest $\mathcal{S}$. Forgetting all the labels on the nodes of $\mathcal{S}$, we get an unlabeled weighted forest of the form $n_1T_1+\cdots +n_rT_r$, where $T_i$ are
mutually distinct unlabeled weighted tree (we call such a tree is of type $\{n_1,\dots, n_r\}$).
Then
$$\Delta_{\mathcal{S}}/G_{\mathcal{S},\underline{\mu}}\cong X^{(n_1)}\times \cdots\times X^{(n_r)}$$
Figure \ref{figure:weighted forest} gives an example of a labeled weighted forest and the corresponding unlabeled weighted forest. The weight $a,b$ are integers.

Therefore we have proved the following decomposition of the Chow motive of $X[n]/\mathfrak{S}_n$:

\begin{thm} \label{motive X[n]/Sn}
For any unordered set of integers $\nu=\{n_1,\dots, n_r\}$ and any integer $m$, let
$\lambda(\nu,m)$ to be the number of unlabeled weighted forest with $n$ leaves, of type $\nu$ and
total weight $m$, such that at each non-leaf $v$ with $c_v$ children, the weight $m_v$ satisfies
$1\le m_v\le (c_v-1)\dim X-1$. Then
$$h(X[n]/\frak{S}_n)=\bigoplus_{\nu,m}\Big{[} h\big{(}X^{(n_1)}\times \cdots\times X^{(n_r)}\big{)}(m)\Big{]}^{\oplus\lambda(\nu,m)}.$$
\end{thm}

\rmk An application of this theorem. MacDonald proved a formula that relates the Betti number of
$X$ and its symmetric powers:
$$\sum_{n=0}^\infty P_tX^{(n)}\cdot T^n=\dfrac{(1+tT)^{b_1}(1+t^3T)^{b_3}\cdots}{(1-T)^{b_0}(1-t^2T)^{b_2}\cdots}$$
where $b_i$ is the $i$-th Betti number of $X$. By the decomposition of the de Rham cohomology of $X[n]/\frak{S}_n$ induced by the motivic
decomposition formula in the above theorem, we can compute the Betti number of $X[n]/\frak{S}_n$ (modulo the combinatorial difficulty of calculating
$\lambda(\nu, m)$).

\egs Here are some examples of $h(X[n]/\mathfrak{S}_n)$ for small $n$. Let $d=\dim X$.

i) n=2. There are $d$ different forests as follows, where each weight $a\in\mathbb{Z}$ ($1\le a\le d-1$) gives a
forest:
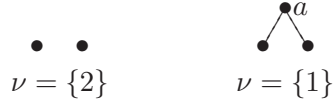
\begin{figure}[h]
\hskip10mm\xy
    (0,0)="m";
    "m"+(0,5)  *{\bullet};
    "m"+(6,5)  *{\bullet};
    "m"+(3,0) *{\nu=\{2\}};
%2
    "m"+(30,0)="m";
    "m"+(3,10)="12"  *{\bullet};
    "m"+(5,10)  *{a};
    "m"+(0,5)="1"  *{\bullet};
    "m"+(6,5)="2"  *{\bullet};
    "12";         "1"**\dir{-};
    "12";         "2"**\dir{-};
    "m"+(3,0) *{\nu=\{1\}};
\endxy
\end{figure}

Therefore $$h(X[2]/\frak{S}_2)\cong h(X^{(2)})\oplus \bigoplus_{a=1}^{d-1}h(X)(a).$$

ii) n=3. The forests are:
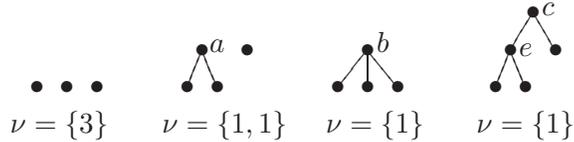
\begin{figure}[h]
\hskip2mm\xy
    (0,0)="m";
    "m"+(0,5)  *{\bullet};
    "m"+(4,5)  *{\bullet};
    "m"+(8,5)  *{\bullet};
    "m"+(3,0) *{\nu=\{3\}};
%2
    "m"+(20,0)="m";
    "m"+(2,10)="12"  *{\bullet};
    "m"+(4,10.5)  *{a};
    "m"+(0,5)="1"  *{\bullet};
    "m"+(4,5)="2"  *{\bullet};
    "m"+(8,10)  *{\bullet};
    "12";         "1"**\dir{-};
    "12";         "2"**\dir{-};
    "m"+(5,0) *{\nu=\{1,1\}};
%3
    "m"+(20,0)="m";
    "m"+(4,10)="123"  *{\bullet};
    "m"+(6,11)  *{b};
    "m"+(0,5)="1"  *{\bullet};
    "m"+(4,5)="2"  *{\bullet};
    "m"+(8,5)="3"  *{\bullet};
    "123";         "1"**\dir{-};
    "123";         "2"**\dir{-};
    "123";         "3"**\dir{-};
    "m"+(5,0) *{\nu=\{1\}};
%4
    "m"+(20,0)="m";
    "m"+(6,15)="123"  *{\bullet};
    "m"+(3,10)="12"  *{\bullet};
    "m"+(9,10)="3"  *{\bullet};
    "m"+(1,5)="1"  *{\bullet};
    "m"+(5,5)="2"  *{\bullet};
    "m"+(8,15.5)  *{c};
    "m"+(5,10)  *{e};
    "123";         "12"**\dir{-};
    "123";         "3"**\dir{-};
    "12";         "1"**\dir{-};
    "12";         "2"**\dir{-};
    "m"+(5,0) *{\nu=\{1\}};
\endxy
\end{figure}

\noindent where the weights $a,b,c,e\in\mathbb{Z}$ satisfy $1\le a,c,e\le d-1$, and $1\le b\le 2d-1$. We have
$$
h(X[3]/\frak{S}_3)\cong h(X^{(3)})\oplus\bigoplus_{i=1}^{d-1}\big{(}h(X^2)(i)\big{)}^{\oplus
3}\oplus\bigoplus_{i=1}^{2d-1}\big{(}h(X)(i)\big{)}^{\oplus min\{i, 2d-i\}}.
$$

iii) n=4. The varieties appear in the decomposition of $h(X[4]/\frak{S}_4)$ are: $$X^{(4)}, X\times
X^{(2)}, X^2, X^{(2)}, X.$$ The decomposition is a bit nasty to be written here. Therefore
we only point out a fact. Consider the forest in Figure \ref{figure:1234}, where $a,b\in\mathbb{Z}$ and $1\le a,b\le d-1$. 
\begin{figure}[h]
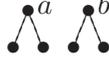

\hskip2mm\xy
    (0,0)="m";
    "m"+(0,5)="1"  *{\bullet};
    "m"+(4,5)="2"  *{\bullet};
    "m"+(8,5)="3"  *{\bullet};
    "m"+(12,5)="4"  *{\bullet};
    "m"+(2,10)="12"  *{\bullet};
    "m"+(10,10)="34"  *{\bullet};
    "m"+(4.2,10.5)  *{a};
    "m"+(12,11)  *{b};
    "12";         "1"**\dir{-};
    "12";         "2"**\dir{-};
    "34";         "3"**\dir{-};
    "34";         "4"**\dir{-};
\endxy
 \caption{\label{figure:1234} An unlabeled weighted forest when $n=4$}
\end{figure}
For any $a<b$, the weighted forest is of type $\nu=\{1,1\}$ and therefore gives a summand $h(X^2)(a+b)$.
However, for $a=b$, this weighted forest has an automorphism exchanging the two trees, thus is of type $\nu=\{2\}$ and  gives a summand
$h(X^{(2)})(2a)$. Due to this kind of automorphism of weighted forests, it seems difficult to compute $\lambda(\nu, m)$.

\begin{question}
 Is there a clean formula for $\lambda(\nu, m)$? (Maybe in terms of a generating function?)
\end{question}

% begin appendix
\appendix

\section{A Formula for the motive of a blow-up} Suppose $f: \Y\to Y$ is the blow-up of a nonsingular projective variety $Y$ along
a nonsingular closed subvariety $V$ of $Y$, and denote by $P$ the exceptional divisor. Denote by $i,j,f,g$ the morphisms as in the following
fibre square
$$\xymatrix{P\ar[r]^j\ar[d]_g\ar@{}[dr]|{\square} & \Y\ar[d]^f\\ V\ar[r]_i & Y.}\label{diamgram_blow-up}$$
Denote by $N:=N_VY$ the normal bundle of $V$ in $Y$. Let $h:=c_1(\O_N(1))\in A^1(P)$. Let $r:=\textrm{codim}_VY$ be the
codimension of $V$ in $Y$.

For $1\le k\le r-1$, define $\alpha_k\in Corr^{-k}(\Y,V)$, $\beta_k\in Corr^{k}(V,\Y)$, $p_k\in Corr^{0}(\Y,\Y)$,
$\alpha_0\in Corr^{0}(\Y,Y)$,  $\beta_0 \in Corr^{0}(Y,\Y)$ and $p_0\in Corr^{0}(\Y,\Y)$ as follows
\begin{equation}\label{many} \left\{
  \begin{array}{ll}
    \alpha_0:=\Gamma_f,\\
    \beta_0:=\Gamma^t_f,\\
     p_0:=\beta_0\circ \alpha_0=\Gamma^t_f\circ\Gamma_f=(f\tm f)^*\Delta_Y,\\
    \alpha_k:=-(j\btm g)_*\Big{(}\sum\limits_{l=0}^{r-1-k} g^*c_{r-1-k-l}(N)h^l\Big{)}\\
    \quad\quad =-(j\btm g)_*\Big{(}\Big{\{}g^*c(N)\dfrac{1}{1-h}\Big{\}}_{r-1-k}\Big{)}, \\
    \beta_k:=(g\btm j)_*h^{k-1}, \\
    p_k:=\beta_k\circ \alpha_k, \\
  \end{array}
\right.\end{equation}
where the subscript $r-1-k$ in the definition of $\alpha_k$
means taking the codimension $(r-1-k)$ component. We will give the proof of the following proposition at the end of this section.

%----------------%
%Main Proposition%
%----------------%
\begin{prop}\label{main1}Define $\alpha_k,\beta_k,p_k,\alpha_0,\beta_0,p_0$  as above. The following holds.

\noindent(i) $\alpha_0\beta_0=\Delta_Y$, $\alpha_k\beta_k=\Delta_V$ for $1\le k\le r-1$; 

$\alpha_i\beta_j=0$ for $i\neq
j$.

\noindent(ii) $p_0, p_1, p_2, ..., p_{r-1}$ are mutually orthogonal projectors of $\Y$, and
$$\sum_{i=0}^{r-1}p_i=\Delta_{\widetilde{Y}} \textrm{ \, in $A(\Y\tm \Y)$,}$$
i.e. equality holds up to rational equivalence.

\noindent(iii) We have the following isomorphisms of motives,
\begin{align*}
&\alpha_0: (\Y, p_0, 0)\simeq h(Y), \textrm{with inverse morphism \,} \beta_0 ,\\
&\alpha_k: (\Y,p_k, 0)\simeq h(V)(k), \textrm{with inverse morphism \,} \beta_k, \textrm{\, for
$1\le k\le r-1$.}
\end{align*}
\end{prop}

Define $\Gamma:=\underset{i=0}{\overset{r-1}{\oplus}}\alpha_i$, $\Gamma':=\underset{i=0}{\overset{r-1}{\sum}}\beta_i$,
then Proposition \ref{main1} can be conveniently reformulated as follows:
\begin{thm}\label{motive for a blow-up}
The correspondence $\Gamma$ gives a canonical isomorphism in $CH\mathcal{M}$,
$$\Gamma: h(\Y)\cong h(Y)\oplus\bigoplus_{k=1}^{r-1}h(V)(k).$$
with an inverse isomorphism given by $\Gamma'.$
\end{thm}

\rmk
When the normal bundle $N$ of $V$ in $Y$ is trivial (for example, when $V$ is a point), $P$ is
isomorphic to a product space $V\tm \P^{r-1}$ and $h=c_1(\O_P(1))$ can be represented (not
canonically) by a product space $H=V\tm \P^{r-2}$ in $P$. In this case, we have simple forms for
the projectors:
\begin{align*}
p_k&=-(j\tm j)_*(H^{r-1-k}\tm_V H^{k-1}), \hbox{ for $1\le k\le r-1$;}\\
p_0&=\Delta+\sum_{k=1}^{r-1}(j\tm j)_*(H^{r-1-k}\tm_V H^{k-1}).
\end{align*}
In general, for a nontrivial
normal bundle $N$, more terms involving the Chern classes of $N$ are needed, and the correspondences cannot
be represented by explicit and natural algebraic cycles.

\rmk\label{compare Manin}
The isomorphism of motives in Theorem \ref{motive for a blow-up} is also a consequence of ``Theorem
on the additive structure of the motif'' of $\Y$ in \cite{Manin} \S9, which states, in our
notation, that there is a split exact sequence
$$\xymatrix{0 \ar[r]& h(V)(r) \ar[r]^-a &h(Y)\oplus h(P)(1)\ar[r]^-b& h(\Y) \ar[r] &0 }.$$
The correspondences appeared in our theorem are not given, at least not explicitly, in Manin's paper.

In order to clarify this point, define
\begin{align*}
&\Phi=c_{r-1}(g^*N/\O_N(-1))\in A^{r-1}(P), c_\Phi=\delta_{P*}(\Phi)\in Corr(P,P),\\
&a=(i_*,c_\Phi\circ g^*),  a'=g_*,\\
&b=f^*+j_*, \hbox{ $b'$ its right inverse,}\\
&d=\Delta_{Y\times P}-aa', d'=\Delta_Y\otimes (\Delta_P-p_0^P) \hbox{ (where
$p_0^P=c_{h^{r-1}}\circ g^*\circ g_*$),}
\end{align*}
denote by $e: \oplus_{k=1}^{r-1}V(k)\to (P,\Delta_P-p_0^P)$  the isomorphism implicitly defined in
\cite{Manin} \S7, and denote by $e'$ the inverse of $e$.

We have the following isomorphisms
\begin{align*}
&\xymatrix{h(Y)\oplus\bigoplus\limits_{k=1}^{r-1}h(V)(k) \ar@<.3ex>[r]^-{\Delta_Y\otimes e}&
  (Y\sqcup P, (\Delta_Y, \Delta_P-p_0^P))\ar@<.3ex>[l]^-{\Delta_Y\otimes e'} \ar@<.3ex>[r]^-d&\ar@<.3ex>[l]^-{d'}}\\
&\xymatrix{(Y\sqcup P, \Delta_{Y\sqcup P}-aa') \ar@<.3ex>[r]^-b &
  (\Y, \Delta_{\Y})\ar@<.3ex>[l]^-{b'}}.
\end{align*}
Hence the following is an isomorphism of Chow motives $$(\Delta_Y\otimes e')\circ d'\circ b':
h(\Y)\cong h(Y)\oplus\bigoplus_{k=1}^{r-1}h(V)\ts \mathbb{L}^k$$ with inverse  $b\circ d\circ
(\Delta_Y\otimes e)$.

Therefore, to write down the correspondence $(\Delta_Y\otimes e')\circ d'\circ b'$, we need to find
explicitly the right inverse $b'$ of $b$. However in \cite{Manin} the construction of $b'$ is based
on the surjectivity of $\gamma: A(\Y\times (Y\sqcup P))\to A(\Y\times \Y)$ as follows:  by the
surjectivity of $\gamma$, there is a cycle class $c\in A(\Y\times (Y\sqcup P))$ (which is not
given, at least explicitly, in \cite{Manin}) such that $\gamma(c)=\Delta_\Y\in A(\Y\times \Y)$.
Then $b'$ is defined to be $(1-aa')c$.
%(it is easy to show that $b'$ does not depend on the choice of $c$).

On the other hand, the correspondences $\Gamma$ and $\Gamma'$ we have constructed in Theorem
\ref{motive for a blow-up} give an explicit construction of $b'$. Indeed,
$b'=d\circ(\Delta_Y\otimes e)\circ\Gamma$.

\begin{proof}[Proof of Proposition \ref{main1}]In the proof, we assume $1\le k\le r-1,\quad 0\le i,j \le r-1.$

The idea is as follows: we study the morphisms
$\alpha_{i*}, \beta_{i*}$ and $p_{i*}$ of Chow groups induced by the
correspondences $\alpha_i, \beta_i$ and $p_i$. As a consequence, the identities of morphisms of
Chow groups which are induced by the identities in Proposition \ref{main1} (i) (ii) hold. On the other hand,
Manin's Identity Principle asserts that the identities of morphisms of Chow groups imply
the identities of correspondences, providing that the correspondences are universal in certain sense. 

By \cite{Voisin} Theorem 9.27, an element $\y\in A(\Y)$ can be expressed uniquely as
$$\y=\sum_{i=1}^{r-1}j_*(g^*a_i\cdot h^{i-1})+f^*y.$$
It is standard to verify

\smallskip

\noindent($\alpha_k$)\quad The morphism $\alpha_{k*}: A(\Y)\to A(V)$ maps
$ \y\mapsto a_{k}.$

\noindent($\beta_k$)\quad The morphism $\beta_{k*}: A(V)\to A(\Y)$ maps
$ x\mapsto j_*(g^*x\cdot h^{k-1}).$

\noindent($\alpha_0$)\quad  The morphism $\alpha_{0*}: A(\Y)\to A(Y)$ maps
$ \y\mapsto y.$

\noindent($\beta_0$)\quad  The morphism $\beta_{0*}: A(Y)\to A(\Y)$ maps
$ y\mapsto f^*y.$

\smallskip

\noindent To give a flavor, we prove only the statement ($\alpha_k$), that is,  $\alpha_{k*}(\y)=a_k$. Define $a_0=-i^*y$. Since $j^*j_*z=-h\cdot z$
for $\forall z\in A(P)$, we have
$$j^*\y=\sum_{i=1}^{r-1}j^*j_*(g^*a_i\cdot h^{i-1})+j^*f^*y=-\sum_{i=1}^{r-1}g^*a_i\cdot h^{i}+g^*i^*y
=-\sum_{i=0}^{r-1}g^*a_i\cdot h^i.$$ By definition (see \cite{Fulton} \S3), the $i$-th Segre class of $N$ is
$$s_{i}(N):=g_*(h^{i+r-1}),$$ hence
\begin{align*}
\alpha_{k*}(\y)&=-g_*\big{(}j^*\y\cdot\sum_{l=0}^{r-1-k} g^*c_{r-1-k-l}(N)\cdot h^l\big{)}\\
&=-g_*\bigg{(}(-\sum_{i=0}^{r-1}g^*a_i\cdot h^i)\cdot\Big{(}\sum_{l=0}^{r-1-k} g^*c_{r-1-k-l}\cdot h^l\Big{)}\bigg{)}\\
&=g_*\Big{(}\sum_{i=0}^{r-1}\sum_{l=0}^{r-1-k}g^*(a_ic_{r-1-k-l})h^{i+l}\Big{)}\\
&=\sum_{i=0}^{r-1}a_i\Big{(}\sum_{l=0}^{r-1-k}c_{r-1-k-l}s_{i+l+1-r}\Big{)}.
\end{align*}
Since $c(N)s(N)=1$ where $c(N):=\sum c_i(N)$ is the total Chern class and $s(N):=\sum s_i(N)$ is the total Segre class, we have
$$\sum_{l=0}^{r-1-k}c_{r-1-k-l}s_{i+l+1-r}=\sum_{l=-\infty}^{+\infty}c_{r-1-k-l}s_{i+l+1-r}=\{c(N)s(N)\}_{i-k}=
\delta_{ik},$$ the first equality is because $s_{i+l+1-r}=0$ for $l<0$ and $c_{r-1-k-l}=0$ for $l>r-1-k$. It follows that $\alpha_{k*}(\y)=a_{k}$, as
we claimed.

\smallskip

The statements ($\alpha_k$),($\beta_k$),($\alpha_0$),($\beta_0$) immediately imply the following identities:
\begin{align*}
&\alpha_{k*}\beta_{k*}=id_{A(V)},\quad \alpha_{0*}\beta_{0*}=id_{A(Y)}, \quad \alpha_{i*}\beta_{j*}=0 \textrm{ for } i\neq j,\\
&(p_ip_j)_*=\delta_{ij}p_{i*}, \quad \sum_{i=0}^{r-1} p_{i*}=id_{A(\Y)}.
\end{align*}

For any smooth scheme $T$, $T\tm \Y$ is the blow-up of $T\tm Y$ along the smooth subvariety $T\tm V$. Denote
$j'=id_T\tm j$, $g'=id_T\tm g$, $f'=id_T\tm f$, $i'=id_T\tm i$, we have the following fiber square:
$$\xymatrix{T\tm P\ar[r]^{j'}\ar[d]^{g'} \ar@{}[dr]|{\square}& T\tm\Y\ar[d]^{f'}\\ T\tm V\ar[r]^{i'} & T\tm Y}$$

We can construct the correspondences $\alpha'_i, \beta'_i, p'_i$ for this fiber square as we did in
(\ref{many}). we have
 $$\alpha'_i=id_T\otimes\alpha_i, \beta'_i=id_T\otimes\beta_i,
p'_i=id_T\otimes p_i.$$

Then (i) and (ii) follows from Manin's Identity Principle.

For (iii), to show that $\alpha_k$ gives an isomorphism $(\Y,p_k,0)\simeq h(V)\ts \mathbb{L}^k$ with
inverse $\beta_k$, we need to show that $p_k=p_k\circ\beta_k\circ\alpha_k$ and
$id=id\circ\alpha_k\circ\beta_k$. but they are direct consequences of the fact that
$\alpha_k\circ\beta_k=\Delta_V$ from (i). The proof for $(\Y, p_0,0)\simeq h(Y)$ is similar.
\end{proof}

\bigskip

\noindent {\sc Li Li }\\
Department of Mathematics\\
University of Illinois at Urbana-Champaign\\
Email: {\tt llpku@math.uiuc.edu}\\

\end{document}